\documentclass[a4paper,12pt]{article} 
\usepackage{amssymb}
\usepackage{amscd}
\usepackage{amsmath} 
\usepackage[dvipdfmx]{graphicx}
\usepackage{latexsym} 
\usepackage{enumerate}

\usepackage{theorem}
\theoremstyle{plain}
\theorembodyfont{\itshape}
\newtheorem{theorem}{Theorem}[section]
\newtheorem{proposition}[theorem]{Proposition}
\newtheorem{lemma}[theorem]{Lemma}
\newtheorem{corollary}[theorem]{Corollary}
\theorembodyfont{\rmfamily}
\newtheorem{definition}[theorem]{Definition}
\newtheorem{remark}[theorem]{Remark}
\newtheorem{example}[theorem]{Example}
\newtheorem{assumption}{Assumption}
\title{\bf Siegel Disks on Rational Surfaces }
\author{Takato Uehara \\ \\
Graduate School of Natural Science and Technology, Okayama University \\
3-1-1 Tsushimanaka, Kita-ku, Okayama 700-8530 
Japan\thanks{E-mail address: {\tt takaue@okayama-u.ac.jp}}} 
\begin{document}
\setlength{\lineskiplimit}{-7pt}
\setlength{\baselineskip}{15.6pt}
\maketitle
\begin{abstract} 
We show the existence of a rational surface automorphism of 
positive entropy with a given number of Siegel disks. 
Moreover, among automorphisms obtained from quadratic birational maps on the 
projective plane fixing irreducible cubic curves, 
we find out an automorphism of positive entropy with multiple Siegel disks. 
\end{abstract} 
\section{Introduction} \label{sec:intro}
A Siegel disk for a holomorphic map on a complex manifold is a domain of the 
manifold preserved by the map such that the restriction to the domain is analytically 
conjugate to an irrational rotation (see Section \ref{sec:pre}). 
Siegel disks are interesting objects and constructed by many authors especially for automorphisms on rational manifolds 
with positive entropy. 
For example, McMullen \cite{M} and Bedford-Kim \cite{BK1, BK2} constructed 
rational surfaces, namely, rational manifolds of dimension $2$, admitting automorphisms of positive entropy with 
Siegel disks by considering a certain class of birational maps on the projective plane. 
Moreover, Oguiso-Perroni \cite{OP} constructed rational manifolds of dimension $\ge 4$ 
admitting automorphisms of positive entropy with an arbitrarily high number of Siegel disks 
by using the product construction made of automorphisms on McMullen's rational surfaces and toric manifolds. 
One of our aims is to show the existence of a rational surface automorphism of positive entropy 
with a given number of Siegel disks. 
\begin{theorem} \label{thm:main1}
For any $k \in \mathbb{Z}_{\ge 0}$, there exists a rational surface $X$ and an automorphism $F : X \to X$ such that 
$F$ has positive entropy $h_{\mathrm{top}}(F)>0$ and $F$ has exactly $k$ fixed points at which Siegel disks are centered. 
\end{theorem}
The automorphism $F$ mentioned in Theorem \ref{thm:main1} is obtained from a birational map 
$f : \mathbb{P}^2 \to \mathbb{P}^2$ of degree $\mathrm{max}\{2,k-1\}$ by blowing up finitely many 
points on the smooth locus of a cubic curve $C$ in $\mathbb{P}^2$. 
When $k \ge 3$, the curve $C$ we considered is the union of three lines meeting a single point. 
\par
Next we consider the case where automorphisms are obtained from quadratic birational maps on 
$\mathbb{P}^2$ that fix a cubic curve $C$. 
Let $f : \mathbb{P}^2 \to \mathbb{P}^2$ be a birational map with its inverse $f^{-1} : \mathbb{P}^2 \to \mathbb{P}^2$ 
and its indeterminacy set $I(f)$, namely, the set of points on which $f$ is not defined. 
We say that $f$ {\it properly fixes} $C$ if 
the indeterminacy sets $I(f^{\pm 1})$ of $f^{\pm 1}$ are both contained in the smooth locus $C^*$ of $C$, 
and $f(C):=\overline{f(C \setminus I(f))}=C$. 
It is known that a certain class of quadratic birational maps properly fixing $C$ is lifted to automorphisms 
with positive entropy by blowing up finitely many points on $C^*$ (see \cite{BK1, BK2, D, M, U1, U2}). 
Let $\mathcal{QF}(C)$ be the set of automorphisms $F : X \to X$ on rational surfaces $X$ 
with positive entropy and 
with the property that there is a quadratic birational map $f : \mathbb{P}^2 \to \mathbb{P}^2$ 
properly fixes $C$ and a blowup $\pi : X \to \mathbb{P}^2$ of points on $C^*$ such that 
the diagonal 
\begin{equation*} 
\begin{CD}
X @> F >> X \\
@V \pi VV @VV \pi V \\
\mathbb{P}^2 @> f >> ~\mathbb{P}^2
\end{CD}
\end{equation*}
commutes. 
In the case where $C$ is non-reduced, 
Bedford-Kim \cite{BK2} constructed $F \in \mathcal{QF}(C)$ with multiple Siegel disks, 
when $C$ is a single line with multiplicity $3$. 
On the other hand, 
McMullen \cite{M} and Bedford-Kim \cite{BK1} constructed $F \in \mathcal{QF}(C)$ with a single Siegel disk, 
when $C$ is reduced but non-irreducible. 
In this article, we focus our attention to the case of irreducible cubic curves, and obtain the following theorem. 
\begin{theorem} \label{thm:main2}
For a reduced irreducible cubic curve $C$ on $\mathbb{P}^2$, if there is an automorphism $F \in \mathcal{QF}(C)$ 
having a Siegel disk, then $C$ is a cuspidal cubic curve. 
Moreover, if $C$ is a cuspidal cubic curve, then $F \in \mathcal{QF}(C)$ admits at most two fixed points at which Siegel disks are centered, 
and there is an automorphism $F \in \mathcal{QF}(C)$ having exactly two fixed points at which Siegel disks are centered. 
\end{theorem}
\par
The existence of a Siegel disk for an automorphism $F$ centered at $x$ implies that 
the derivative $DF(x)$ of $F$ at $x$ has 
multiplicatively independent eigenvalues $(\mu,\nu)$ with $|\mu|=|\nu|=1$ (see Section \ref{sec:pre}). 
Conversely, results from transcendence theory guarantee that $F$ has a Siegel disk centered at $x$ 
under the assumption that the multiplicatively independent eigenvalues 
$(\mu,\nu)$ with $|\mu|=|\nu|=1$ are algebraic. 
Moreover, if algebraic eigenvalues $(\mu,\nu)$ with $|\mu|=|\nu|=1$ have Galois conjugates 
$(\mu_*,\nu_*)$ satisfying $|\mu_* \nu_*|=1$ but $|\mu_*/ \nu_*| \neq 1$, 
then $(\mu,\nu)$ are multiplicatively independent (see also \cite{M}). 
Our task is thus to construct automorphisms whose derivatives have such a pair 
$(\mu,\nu)$ of eigenvalues. 
Note that in our construction, the automorphisms are obtained from birational maps, 
and the birational maps considered here have explicit forms with parameters. 
\par
After preliminary studies in Section \ref{sec:pre}, 
Sections \ref{sec:Quad} and \ref{sec:lines} are devoted to constructing automorphisms with Siegel disks 
in order to prove Theorems \ref{thm:main2} and \ref{thm:main1} respectively, 
and Sections \ref{sec:pp1} and \ref{sec:pp2} are devoted to proving two propositions needed in our discussion. 
\\[2mm]
{\bf Acknowledgment}. 
The author is supported by Grant-in-Aid for Young Scientists (B) 24740096. 
\section{Preliminary} \label{sec:pre}
In this section, we briefly review some well-known facts about 
Siegel disks on complex surfaces, automorphisms on rational surfaces and cubic curves on the projective plane used later. 
We refer to \cite{D, M, U1, U2}, in which many of the results are proved. 
\par 
First we recall a Siegel disk on a complex surface (see \cite{M}). 
For a unit disk $\Delta^2:=\{(x,y) \in \mathbb{C}^2 \, | \, |x| \le 1, |y| \le 1\}$, 
a linear automorphism $L : \Delta^2 \to \Delta^2$ given by $L(x,y)=(\mu x, \nu y)$ is called 
an {\it irrational rotation} if $|\mu|=|\nu|=1$ and $(\mu,\nu)$ are {\it multiplicatively independent}, that is, 
they satisfy $\mu^k \nu^{l} \neq1$ for any $(k,l) \neq (0,0) \in \mathbb{Z}^2$. 
\begin{definition} \label{def:Siegel}
Let $X$ be a complex surface and $F$ be an automorphism on $X$. 
A domain $U \subset X$ is called a {\it Siegel disk} for $F$ centered at $p \in U$ 
if $F(U)=U$ and $F : (p,U) \to (p,U)$ is analytically conjugate 
to an irrational rotation $L : (0,\Delta^2) \to (0,\Delta^2)$.  
\end{definition}
It is obvious that the derivative $DF(p)$ of $F$ at $p$ is an irrational rotation when $F$ has a Siegel disk centered at $p$. 
Conversely, results from transcendence theory say that if $DF(p)$ is an irrational rotation with {\it algebraic eigenvalues}, 
then $F$ has a Siegel disk centered at $p$ (see \cite{M}). 
\par 
Next we consider rational surfaces. Here we assume that a rational surface $X$ admits a birational morphism $\pi : X \to \mathbb{P}^2$ 
(see \cite{D, M, U1, U2}). 
Then it is known that $\pi$ is expressed as a composition 
\[
\pi : X=X_\rho \overset{\pi_{\rho}}{\longrightarrow} X_{\rho-1} 
\overset{\pi_{\rho-1}}{\longrightarrow} \cdots \overset{\pi_{2}}{\longrightarrow}
X_{1} \overset{\pi_1}{\longrightarrow} X_0=X, 
\]
where $\pi_i : X_{i} \to X_{i-1}$ is the blowup of a point $p_i \in X_{i-1}$ with the exceptional curve 
$\mathcal{E}_i :=\pi_i^{-1}(\{p_i\})$, which is isomorphic to $\mathbb{P}^1$. 
Since $\pi_i$ induces an isomorphism 
$\pi_i|_{X_i \setminus \mathcal{E}_i} : X_i \setminus \mathcal{E}_i \to X_{i-1} \setminus \{p_{i}\}$, 
we will identify each point $x \in X_i \setminus \mathcal{E}_i$ with 
$\pi_{i}(x) \in X_{i-1} \setminus \{p_{i}\}$ in this article. 
Moreover, if $p$ is a point on an exceptional curve,  
we sometimes says that $p$ is an {\it infinitely near point} on $\mathbb{P}^2$, 
or a {\it point} on $\mathbb{P}^2$ for short. 
On the other hand, a point is said to be {\it proper} if it is {\it not} an infinitely near point.  
The total transform $E_i := \pi_{\rho}^* \circ \cdots \circ \pi_{i+1}^*(\mathcal{E}_i)$ is called 
the {\it exceptional divisor} over $p_i$. 
Then $\pi$ gives an expression of the cohomology group : 
\[
H^2(X;\mathbb{Z}) \cong \mathrm{Pic}(X)=\mathbb{Z} [H] \oplus \mathbb{Z} [E_1] 
\oplus \cdots \oplus \mathbb{Z} [E_\rho], 
\]
where $H$ is the total transform $\pi^*(L)$ of a line $L$ in $\mathbb{P}^2$.  
The intersection form on the cohomology group $H^2(X;\mathbb{Z})$ is given by 
\[
\left\{
\begin{array}{ll}
([H],[H])=1 & ~ \\[2mm]
([E_i],[E_j])=-\delta_{i,j} & (i,j=1,\dots, \rho) \\[2mm]
([H],[E_i])=0 & (i=1,\dots, \rho). 
\end{array}
\right. 
\]
\par
Let $F : X \to X$ be an automorphism on $X$. 
Then $F$ induces the action $F^* : H^2(X;\mathbb{Z}) \to H^2(X;\mathbb{Z})$ on the cohomology group. 
By theorems of Gromov and Yomdin, the {\it topological entropy of $F$} is given by 
$h_{\mathrm{top}}(F)= \log \lambda(F^*) \ge 0$ with the spectral radius $\lambda(F^*)$ of $F^*$. 
Moreover since $F^*$ preserves the K\"ahler cone and the intersection form with signature $(1,\rho)$, 
it is seen that the characteristic polynomial of $F$ is 
expressed as 
\[
\mathrm{det} (t I - F^*) = 
\left\{
\begin{array}{ll}
R_{F}(t) \quad & (\lambda(F^*)=1) 
\\[2mm]
R_{F}(t) S_{F}(t) \quad & (\lambda(F^*)>1), 
\end{array}
\right.
\]
where $R_{F}(t)$ is a product of cyclotomic polynomials, and 
$S_{F}(t)$ is a Salem polynomial, namely, the minimal polynomial 
of a Salem number. 
Here, a {\it Salem number} is an algebraic unit $\delta>1$ 
such that its conjugates include $\delta^{-1} < 1$ 
and the conjugates other than $\delta^{\pm 1}$ lie on the unit circle. 
Hence if $\lambda(F^*)>1$ then it is a root of $S_{F}(t)=0$. 
\par 
Now we consider a cubic curve $C \subset \mathbb{P}^2$, that is, 
a reduced (possibly non-irreducible or singular) curve of degree three, 
with its smooth locus $C^*$ (see \cite{D, M}). 
Denote by $\mathrm{Pic}^0(C) \subset \mathrm{Pic}(C)$ the subgroup consisting of 
divisor classes whose restrictions to each irreducible component of $C$ have degree zero. 
Then it is known that $\mathrm{Pic}^0(C) \cong \mathbb{C} / \Gamma$, where $\Gamma \subset \mathbb{C}$ 
is a lattice with rank given by either 
\begin{enumerate}
\item $\mathrm{rank}~\Gamma = 2$ if $C$ is smooth; or 
\item $\mathrm{rank}~\Gamma = 1$ if $C$ is a nodal cubic, or a conic with a transverse line, or three lines meeting in three points; or 
\item $\mathrm{rank}~\Gamma = 0$ if $C$ is a cuspidal cubic, or a conic with a tangent line, or three lines through a single point.  
\end{enumerate}
Let $V_1,\dots,V_r$ be the irreducible components of $C$. Note that $1 \le r \le 3$ as $C$ is a cubic curve. 
Moreover fix points $0_i \in V_i \cap C^*$ so that $\sum_{i=1}^r \mathrm{deg} V_i \cdot [0_i]=0$, namely, 
the divisor $\sum_{i=1}^r \mathrm{deg} V_i \cdot 0_i$ is 
the restriction of a line $L \subset \mathbb{P}^2$ to $C^*$, 
where $\mathrm{deg} V_i \in \mathbb{Z}_{>0}$ is the degree of the component $V_i$ in $\mathbb{P}^2$. 
For each $1 \le j \le r$, let $\kappa : V_j \cap C^* \to \mathrm{Pic}^0(C)$ be the map defined by $\kappa(p)=[p]-[0_j]$. 
Then $\kappa$ is a bijection, 
which gives the group structure on $V_j \cap C^*$ isomorphic to $\mathrm{Pic}^0(C) \cong \mathbb{C} / \Gamma$, 
with the property that three points $q_1,q_2,q_3 \in C^*$ satisfy $\sum_{i=1}^3 [q_i]= 0$ if and only if 
$\sum_{i=1}^3 \kappa(q_i)=0$ and $\# \{ i \, | \, q_i \in V_j \}=\mathrm{deg} V_j$ for any $1 \le j \le r$ (see \cite{D}). 
\par 
Let $f : \mathbb{P}^2 \to \mathbb{P}^2$ be a birational map on $\mathbb{P}^2$. 
In general, $f$ admits the {\it indeterminacy set} $I(f)$, namely, the finite set on which $f$ cannot be defined (see \cite{U1, U2}). 
Note that $I(f)$ is a cluster, that is, if $p \in I(f)$ is infinitely near to a point $q$ then $q \in I(f)$. 
All birational maps considered in this article are assumed to belong to the set $\mathcal{B}(C)$ 
of birational maps $f$ {\it properly fixing} $C$, namely, $I(f^{\pm 1}) \subset C^*$ and $f(C)=C$. 
Here, if $I(f^{\pm 1})$ contain an infinitely near point $p$, then $p \in C^*$ means that 
$p$ belongs to the strict transform $\pi^{-1}(C^*)$, where $\pi : X \to \mathbb{P}^2$ is a birational morphism 
such that $p$ is proper on $X$. 
When $f \in \mathcal{B}(C)$, there is $\delta(f) \in \mathbb{C}^*$, called the {\it determinant} of $f$, 
such that $f^* \eta=\delta (f) \eta$, where $\eta$ is a nowhere vanishing meromorphic $2$-form on $\mathbb{P}^2$ 
with a simple pole along $C$. 
The determinant $\delta(f)$ satisfies $\delta(f)=\mathrm{Det} Df(p)$ for 
any fixed point $p \in \mathbb{P}^2 \setminus C$ of $f$. 
Moreover it should be noted that $f$ preserves the smooth locus $C^*$ under our assumption. 
Thus $f$ induces the actions $f_* : \mathrm{Pic}(C) \to \mathrm{Pic}(C)$ 
and $f_* : \mathrm{Pic}^0(C) \to \mathrm{Pic}^0(C)$. 
Through the Poincar\'e residue map, it turns out that the action $f_*$ on 
$\mathrm{Pic}^0(C) \cong \mathbb{C} / \Gamma$ is given by $f_*(t)= \delta(f) t$ for $t \in \mathbb{C} / \Gamma$ (see \cite{M}). 
Note that if $\mathrm{rank}~\Gamma \ge 1$ then $\delta(f)$ must be a root of unity as $\delta(f) \Gamma=\Gamma$, 
while if $\mathrm{rank}~\Gamma =0$ then $\delta(f)$ may be an arbitrary nonzero complex number. 
\par 
One of our interests is to construct automorphisms on rational surfaces. 
From birational maps on $\mathbb{P}^2$ satisfying a certain assumption, 
we obtain rational surface automorphisms. 
\begin{proposition} \label{prop:auto}
Assume that $C \subset \mathbb{P}^2$ is a reduced cubic curve. 
\begin{enumerate} 
\item For a birational map $f : \mathbb{P}^2 \to \mathbb{P}^2$ in $\mathcal{B}(C)$, 
assume that any indeterminacy point $p \in I(f^{-1})$ satisfies $f^{m}(p) \in I(f)$ for some $m=m(p) \ge 0$. 
Then there is a blowup $\pi : X \to \mathbb{P}^2$ of points on $C^*$ such that 
$\pi$ lifts $f : \mathbb{P}^2 \to \mathbb{P}^2$ to an automorphism $F : X \to X$. 
\item Assume that a birational map $f : \mathbb{P}^2 \to \mathbb{P}^2$ in $\mathcal{B}(C)$ 
is lifted to an automorphism $F : X \to X$ by a blowup $\pi : X \to \mathbb{P}^2$ of 
points on $C^*$. 
Then any indeterminacy point $p \in I(f^{-1})$ satisfies $f^{k}(p) \notin I(f)$ with $0 \le k < m_p$ 
and $f^{m_p}(p) \in I(f)$ for some $m_p \ge 0$. 
Moreover, $\pi$ admits an expression $\pi=\pi_0 \circ \nu : X \to \mathbb{P}^2$, 
where $\pi_0 : X_0 \to \mathbb{P}^2$ is the blowup of the points 
$\{f^k(p) \, | \, p \in I(f^{-1}), 0 \le k \le m_p \}$ on $C^*$, that lifts $f$ to an automorphism $F_0 : X_0 \to X_0$, 
and $\nu : X \to X_0$ is a birational morphism.   
\end{enumerate}
\end{proposition}
{\it Proof}. 
(1) (see \cite{U1}). Let $(p, q) \in I(f^{-1}) \times I(f)$ be a pair of {\it proper} points such that 
$f^n(p)=q$ with $n = \min \{ m \in \mathbb{N} \, | \, f^m(p') =q' \text{ for }(p',q') \in I(f^{-1}) \times I(f)\}$. 
Under our assumption, such a pair $(p,q)$ exists, and 
from the minimality of $n$, the orbit $\{f^i(p)\}_{i = 0}^n$ consists of distinct proper points 
on the smooth locus $C^*$. 
Now let $X_0 \to \mathbb{P}^2$ be the blowup of $\{f^i(p)\}_{i = 0}^n$. 
The blowup lifts $f : \mathbb{P}^2 \to \mathbb{P}^2$ to a birational map $f_0 : X_0 \to X_0$, 
which satisfies  
\[
I(f_0^{-1}) = I(f^{-1}) \setminus \{p\}, \quad I(f_0) = I(f) \setminus \{q\}. 
\]
Note that $\# I(f^{-1})=\# I(f)$. 
Therefore as long as $\# I(f_0^{-1})=\# I(f_0) >0$, one can repeat the argument 
by replacing $f : \mathbb{P}^2 \to \mathbb{P}^2$ with $f_0 : X_0 \to X_0$. 
In the end, a resulting map becomes an automorphism. 
See \cite{U1} for a more detailed discussion. \\[2mm]
(2) (see \cite{U2}). 
We notice that if $p \in I(f^{-1})$ satisfies $f^{k}(p) \notin I(f)$ for $0 \le k \le m-1$ then 
$f^m(p)$ is a well-defined point in $I(f^{-m})$. 
As $\pi$ lifts $f^m$ to the automorphism $F^m$, the point $f^m(p)$ must be blown up by $\pi$. 
Since the number of points blown up by $\pi$ is finite, there is $m_p \ge 0$ such that 
$f^k(p) \notin I(f)$ for $0 \le k \le m_p-1$ and $f^{m_p}(p) \in I(f)$. 
Moreover, $\pi$ blows up the points $\{f^k(p)\, | \, p \in I(f^{-1}), 0 \le k \le m_p \}$, 
and hence $\pi$ admits the expression $\pi=\pi_0 \circ \nu : X \to \mathbb{P}^2$. 
The blowup $\pi_0$ lifts $f$ to an automorphism $F_0$ from a similar argument in the proof of (1). 
See \cite{U2} for a more detailed discussion. 
\hfill $\Box$ \par\medskip 
\begin{definition} \label{def:proper}
For a birational map $f \in \mathcal{B}(C)$ satisfying the assumption in Proposition \ref{prop:auto} (1), 
the blowup $\pi_0$ given in Proposition \ref{prop:auto} (2) is called the {\it proper blowup} for $f$. 
\end{definition}
\begin{remark} \label{rem:iden}
Let $f : \mathbb{P}^2 \to \mathbb{P}^2$ be a birational map lifted to an automorphism $F : X \to X$ 
by a blowup $\pi : X \to \mathbb{P}^2$. 
With the identification of a point $p \in X$ with $\pi(p) \in \mathbb{P}^2$ under the assumption that $\pi(p) \notin I(\pi^{-1})$, 
the dynamical behaviour of $F$ around $p$ is the same as that of $f$ around the corresponding point. 
In particular, $F$ has a Siegel disk centered at $p$ if and only if so does $f$. 
\end{remark}
\par 
The next lemma is used to calculate the cohomological actions of automorphisms. 
\begin{lemma} \label{lem:exchange}
Let $\pi$ be the proper blowup for $f$, that lifts $f$ to an automorphism $F$,and let $p_1, \dots, p_{\rho}$ be the points blown up by $\pi$ 
and $E_l$ be the exceptional divisor over $p_l$. 
If a point $p_i$ satisfies $p_i \notin I(f^{-1})$, then the action $F^*$ of $F$ sends $E_i$ to $E_{j}$ for some $j \neq i$.  
\end{lemma}
{\it Proof}. 
Under the notations given in the proof of Proposition \ref{prop:auto} (1), 
we may assume that $p_i=f^k(p) \notin I(f^{-1})$ for some $k \ge 1$, 
as the other cases can be treated in a similar manner. 
Note that $f^m(p) \notin I(f^{-1})$ for any $0 \le m \le k$ in this case. 
As is mentioned in the proof of Proposition \ref{prop:auto} (1), 
the blowup $X_0 \to \mathbb{P}^2$ of $\{f^i(p)\}_{i = 0}^n$ lifts $f$ to $f_0 : X_0 \to X_0$, and then $f_0$ sends 
$\mathcal{E}^{k-1}$ to $\mathcal{E}^{k}$, where $\mathcal{E}^{l}$ is the exceptional curve over $f^{l}(p)$. 
As the indeterminacy set is a cluster, any point on $\mathcal{E}^{k}$ is not an indeterminacy point of $f^{-1}$. 
Moreover, since $\pi$ is a proper blowup for $f$, there is a point $p' \in \mathcal{E}^k$ blown up by $\pi$ 
if and only if there is a point $p'' \in \mathcal{E}^{k-1}$ blown up by $\pi$ such that $f_0(p'')=p'$, 
which shows that $F$ sends the irreducible components of the exceptional divisor over $f^{k-1}(p)$ to those over $f^{k}(p)$. 
Therefore $F^*$ sends the exceptional divisor over $f^k(p)$ to that over $f^{k-1}(p)$. 
\hfill $\Box$ \par\medskip 
\begin{example} \label{ex:quad}
We consider a quadratic birational map on $\mathbb{P}^2$. 
It is known that the inverse of any quadratic birational map is also quadratic, 
and the indeterminacy set of a quadratic birational map consists of exactly three non-collinear
(possibly infinitely near) points. 
Let $f : \mathbb{P}^2 \to \mathbb{P}^2$ be a quadratic birational map in $\mathcal{B}(C)$, 
and put $I(f^{\pm 1})=\{p_1^{\pm}, p_2^{\pm},p_3^{\pm}\} \subset C^*$. 
Then $f$ lifts to an automorphism if and only if $f^{k}(p_i^{-}) \notin I(f)$ for $0 \le k < n_i$ and 
$f^{n_i}(p_i^{-})=p_{\sigma(i)}^+$ for any $i \in \{1,2,3\}$, 
where $n_1, n_2,n_3 \ge 0$ are integers and $\sigma : \{1,2,3\} \to \{1,2,3\}$ is a permutation. 
Let $\pi_0$ be the proper blowup for $f$, which lifts $f$ to an automorphism $F_0 : X_0 \to X_0$. 
With a suitable matching of the indices between forward and backward indeterminacies, 
the action $F_0^* : H^2(X_0 ; \mathbb{Z}) \to H^2(X_0 ; \mathbb{Z})$ is expressed as 
\[
\left\{
\begin{array}{cll}
[H]  &\mapsto  2[H] - [E_1^{n_1}] - [E_2^{n_2}] - [E_3^{n_3}] & ~ \\[2mm]
[E_i^0]  & \mapsto  [H] - [E_{\sigma(j)}^{n_j}] - [E_{\sigma(k)}^{n_k}] & (\{ i,j,k \}= \{ 1,2,3 \}) \\[2mm]
[E_l^{m}] & \mapsto  [E_l^{m-1}] & (l \in \{1,2,3\}, m \ge 1), 
\end{array}
\right.
\]
where $E_{l}^m$ is the exceptional divisor over $f^m(p_l^-)$ (see \cite{D, U1}). 
\end{example}
\par 
As is mentioned in Proposition \ref{prop:auto}, we assume that the points $(p_1, \dots, p_\rho)$ 
blown up by $\pi : X \to \mathbb{P}^2$ lie on the smooth locus $C^*$ of the cubic curve $C$, 
and also assume that $\pi$ lifts a birational map $f : \mathbb{P}^2 \to \mathbb{P}^2$ in 
$\mathcal{B}(C)$ to an automorphism $F : X \to X$.
Since $f$ preserves $C$, the automorphism $F$ also preserves 
the strict transform $Y:= \pi^{-1}(C)$ of $C$. 
Moreover, as the points $p_i$ lie on $C^*$, the curve $Y$ is isomorphic to $C$ and anticanonical on $X$, namely, $[Y]=-K_X$, where 
$K_X:=-3 [H] + \sum_{i=1}^\rho [E_i]$. 
Under the above notation, we have the following proposition. 
\begin{proposition} \label{prop:eigen}
Assume that $\mathrm{Pic}^0(C) \cong \mathbb{C}$, and also assume that
\begin{enumerate}
\item $\# \{ 1 \le i \le \rho \, | \, p_i \in V_j\} \ge \mathrm{deg} V_j$ for any irreducible component 
$V_j$ of $C$, and 
\item $\kappa(p_i) \neq 0$ for some $1 \le i \le \rho$, where $\kappa : V_j \cap C^* \to \mathrm{Pic}^0(C) \cong \mathbb{C}$ 
is given by $\kappa(p)=[p]-[0_j]$. 
\end{enumerate}
Then, the determinant $\delta(f)$ is an eigenvalue of $F^* : H^2(X;\mathbb{Z}) \to H^2(X;\mathbb{Z})$. 
\end{proposition}
{\it Proof}. 
Let $r \in \{1,2,3\}$ be the number of irreducible components of $C$. 
From the assumption (1) we may assume that $\# \{ 1 \le i \le 3 \, | \, p_i \in V_j \} = \mathrm{deg} V_j$ for $1\le j \le r$, 
after reordering $(p_i)$ if necessary, and also choose  
$\sigma : \{1,\dots,\rho\} \to \{1,\dots,r\}$ so that $p_i \in V_{\sigma(i)}$ for $1 \le i \le \rho$. 
Let us consider the restriction map 
$u : H^2(X;\mathbb{Z}) \cong \mathrm{Pic}(X) \to \mathrm{Pic}(Y) \cong \mathrm{Pic}(C)$, explicitly given by 
\[
u[H]= \sum_{i=1}^r \mathrm{deg}V_i \cdot [0_i], \qquad u[E_i]=[p_i] \quad (i=1,\dots, \rho). 
\]
Then the following diagram commutes: 
\[
\begin{CD}
H^2(X;\mathbb{Z}) @> F_* >> H^2(X;\mathbb{Z}) \\
@V u VV @VV u V \\
\mathrm{Pic}(C) @> f_* >> ~\mathrm{Pic}(C). 
\end{CD}
\]
For simplicity, we denote by the same notation $V_i$ the strict transform $\pi^{-1}(V_i)$. 
Since $F^*$ preserves the intersection form and permutes the curves $\{V_1,\dots,V_r\}$, 
it preserves the orthogonal complement $\mathcal{H}_X := \{[V_1],\dots,[V_r] \}^{\perp} \subset H^2(X;\mathbb{Z})$, 
generated by $(B_0,B_{r+1},\dots, B_{\rho})$ with 
\[
B_0 := [H] - [E_1] - [E_2] - [E_3], \qquad 
B_i := [E_i] - [E_{\sigma(i)}], \quad (i=r+1,\dots,\rho).  
\]
We notice that the image of $u$ restricted to $\mathcal{H}_X$ is contained in $\mathrm{Pic}^0(C)$.  
\par 
Now let us fix a vector $\xi \in H^2(X;\mathbb{C})= H^2(X;\mathbb{Z}) \otimes \mathbb{C}$ satisfying 
\[
\kappa(p_i)= - (\xi, [H]/3-[E_i]) \in \mathrm{Pic}^0(C) \cong \mathbb{C}. 
\]
Note that under the assumption (2), the vector $\xi$ is nonzero and unique in $H^2(X;\mathbb{C})/ \mathbb{C} [K_X]$. 
Then we have 
\[
u(B_0) = \sum_{i=1}^r \mathrm{deg} V_i \cdot [0_i] - \sum_{i=1}^3 [p_i] = \sum_{i=1}^3 \{ [0_{\sigma(i)}] - [p_i]  \} 
= - \sum_{i=1}^3 \kappa(p_i) = \sum_{i=1}^3 (\xi, [H]/3-[E_i]) = (\xi, B_0). 
\]
In a similar manner, it follows that $u(B_i)=(\xi,B_i)$ and thus $u(D)=(\xi,D)$ for any $D \in \mathcal{H}_X$. 
Note that the action $f_*$ on $\mathrm{Pic}^0(C) \cong \mathbb{C}$ is given by $f_*(t)= \delta(f) t$ for $t \in \mathbb{C}$. 
Therefore for any $D \in \mathcal{H}_X$, we have 
\[
\begin{array}{l} 
u(F_* D) = (\xi,F_* D) =(F^* \xi,D) \\
~~~ = f_* u(D) = \delta (f) (\xi,D) = (\delta (f) \xi,D), 
\end{array}
\]
which yields $F^* \xi = \delta(f) \xi + \sum_{i=1}^r c_i [V_i]$ for some $c_i \in \mathbb{C}$. 
Since $F^*$ preserves $\{[V_1],\dots, [V_r] \}$, 
$\delta(f)$ is an eigenvalue of $F^*$. 
The proposition is established. 
\hfill $\Box$ \par\medskip 
Now, in addition to the assumptions in Proposition \ref{prop:eigen}, we also assume that
$C$ is a cuspidal cubic curve and the determinant $\delta(f)$ is not a root of unity. 
Then $\delta(f)$ is a root of the Salem polynomial $S_F(t)=0$ by Proposition \ref{prop:eigen}, 
and the entropy of $F$ is positive: $h_{\mathrm{top}}(F)= \log \lambda(F^*) > 0$. 
In this case, the birational morphism $\nu: X \to X_0$ mentioned in Proposition \ref{prop:auto} is expressed as follows. 
Let $q \in Y^*$ be a fixed point on the smooth locus $Y^* \cong \mathbb{C}$ of the anticanonical curve $Y$, 
which uniquely exists as $F$ has the determinant $\delta(f) \neq 1$. 
A result in \cite{U2} says that if $\nu$ is not an isomorphism, then there is a unique 
$(-1)$-curve passing through $q$, which is contracted by $\nu$ and is preserved by $F$. 
Through the contraction of the $(-1)$-curve, $F$ descends to an automorphism. 
Repeating this argument,  we can consider the decomposition 
\begin{equation} \label{eqn:comon}
\nu : X=X_m \overset{\nu_m}{\longrightarrow} X_{m-1} 
\overset{\nu_{m-1}}{\longrightarrow} \cdots \overset{\nu_{2}}{\longrightarrow}
X_{1} \overset{\nu_1}{\longrightarrow} X_0, 
\end{equation}
where $\nu_i : X_{i} \to X_{i-1}$ is the contraction of a $(-1)$-curve 
through $p_{i}$ to $p_{i-1}$ with $p_m:=q$. 
Then $F$ descends to an automorphism $F_0 : X_0 \to X_0$. 
\par
Let $\mathcal{N}_i \subset X$ be the strict transform of 
the exceptional curve of $\nu_i$ under $\nu_{i+1} \circ \cdots \circ \nu_m$. 
As $\mathcal{N}_i$ is isomorphic to $\mathbb{P}^1$ and is preserved by $F$, 
we inductively let $q_i$ be the unique fixed point on $\mathcal{N}_i \setminus \{ q_{i+1} \}$ of $F$ with $q_{m+1}:=q$. 
In particular, $(q_1,\dots,q_{m}, q)$ are all of the fixed points lying on the exceptional divisors of $\nu$. 
Moreover, let $p \in C$ be the singular point of $C$, which is also a fixed point of $F$. 
\begin{proposition}[\cite{U2}] \label{pro:eigen}
Under the above assumptions, we have the following. 
\begin{enumerate}
\item The eigenvalues of $DF$ at $p$ 
are $1/\delta(f)^2$ and $1/\delta(f)^{3}$. 
\item The eigenvalues of $DF$ at $q$ 
are $\delta(f)$ and $1/\delta(f)^{N-4}$, 
where $N= \mathrm{rank}~\mathrm{Pic}(X)$. 
\item The eigenvalues of $DF$ at $q_i$ for $1 \le i \le m$ 
are $\delta(f)^{N-m+i-4}$ and $1/\delta(f)^{N-m+i-5}$. 
\end{enumerate}
In particular, $F$ has no Siegel disk centered at any fixed point on the anticanonical curve $Y$ and 
the exceptional divisors of $\nu$.  
\end{proposition}
\par 
Next we give an estimate of the number of isolated fixed points of an automorphism. 
\begin{proposition} \label{prop:fpf}
Assume that an automorphism $F : X \to X$ on a rational surface $X$ has positive entropy, and 
the derivative $DF(x)$ of $F$ on any fixed point $x$ has an eigenvalue different from $1$. 
Then $F$ has at most $\mathrm{Tr} (F^*|_{H^2(X;\mathbb{Z})} )+2$ isolated fixed points. 
\end{proposition}
We postpone its proof to Section \ref{sec:pp1}. 
The following two propositions are applications of Proposition \ref{prop:fpf}. 
\begin{proposition} \label{pro:nof}
Let $C \subset \mathbb{P}^2$ be a reduced cubic curve with $\mathrm{Pic}^0(C) \cong \mathbb{C}$, and 
let $F : X \to X$ be an automorphism with positive entropy such that $F$ is obtained from 
a birational map $f \in \mathcal{B}(C)$ by the blowup $\pi : X \to \mathbb{P}^2$ of points on $C^*$. 
Assume that $\delta(f)$ is not a root of unity. 
Then $F$ has at most $\mathrm{Tr} (F^*|_{H^2(X;\mathbb{Z})} )+2$ isolated fixed points. 
\end{proposition}
{\it Proof}. 
First we notice that our assumption says that for any fixed point $x$ on 
the anticanonical curve $Y=\pi^{-1}(C)$, the derivative 
$DF(x)$ of $F$ on $x$ has an eigenvalue different from $1$. 
Indeed, if $x$ lies on the smooth locus $Y^*$, then $DF(x)$ has $\delta(f)$ as an eigenvalue. 
On the other hand, if $x$ is a singular point of $Y$, then $DF(x)$ has eigenvalues of the form $\epsilon \delta(f)^{-m}$, 
where $\epsilon$ is a root of unity and $m \in \mathbb{Z}_{>0}$ is a positive integer (see \cite{M}, \S 9).  
\par
This remains true for any fixed point $x$ outside $Y$, 
since $\mathrm{Det} DF(x)=\delta(f) \neq 1$ from the existence of a nowhere vanishing meromorphic $2$-form 
$\eta_X= \pi^* \eta$ on $X$ with $(\eta_X)= -Y$ and $F^* \eta_X=\delta(f) \eta_X$. 
Hence the proposition follows from Proposition \ref{prop:fpf}. 
\hfill $\Box$ \par\medskip 
\begin{proposition} \label{pro:sieg}
For a cuspidal cubic curve $C$, let $f \in \mathcal{B}(C)$ be a quadratic birational map 
with $\delta(f)$ being not a root of unity such that $f$ is lifted to an automorphism $F : X \to X$ 
by the blowup $\pi : X \to \mathbb{P}^2$  of points on $C^*$. 
Then $F$ has at most $2$ fixed points at which Siegel disks are centered. 
\end{proposition}
{\it Proof}. 
Note that $\pi$ satisfies the assumptions in Proposition \ref{prop:eigen}. 
Indeed the assumption (1) holds as it follows form Proposition \ref{prop:auto} (2) that three indeterminacy points $\{p_1^+, p_2^+,p_3^+\}$ of $f$ 
are blown up by $\pi$. 
Moreover the assumption (2) also holds as the points $\{p_1^+, p_2^+,p_3^+\}$ are not collinear. 
Hence Proposition \ref{prop:auto} (2) and the above argument show 
that the blowup $\pi$ can be decomposed as $\pi= \pi_0 \circ \nu$, where 
$\pi_0 : X_0 \to \mathbb{P}^2$ is the proper blowup for $f$, which lifts $f$ to an automorphism 
$F_0 : X_0 \to X_0$,  and $\nu : X \to X_0$ is expressed as the decomposition (\ref{eqn:comon}). 
The cohomological action $F_0^* : H^2(X_0;\mathbb{Z}) \to H^2(X_0;\mathbb{Z})$ is given in 
Example \ref{ex:quad}, which means that $\mathrm{Tr}(F_0^*|_{H^2(X_0;\mathbb{Z})}) \le 2$. 
Hence $F_0$ has at most $4$ isolated fixed points by Proposition \ref{pro:nof}, since $h_{\mathrm{top}}(F_0)= h_{\mathrm{top}}(F) > 0$. 
Among the fixed points, two fixed points lie on the anticanonical curve $Y_0=\pi_0^{-1}(C)$ of $X_0$, 
at which no Siegel disks are centered from Proposition \ref{pro:eigen}. 
On the other hand, Proposition \ref{pro:eigen} also shows that at none of 
the fixed points of $F$ on the exceptional divisors of $\nu$, a Siegel disk is centered. 
Since each fixed point of $F$ either is identified with that of $F_0$ or lies on the exceptional divisors of $\nu$ (see also Remark \ref{rem:iden}), 
$F$ has at most $2$ fixed points at which Siegel disks are centered. 
\hfill $\Box$ \par\medskip 
We conclude this section by stating a result for a class of birational maps with algebraic coefficients 
that we will treat in the following sections.
To this end, for a reduced cubic curve $C \subset \mathbb{P}^2$ and a birational map $f : \mathbb{P}^2 \to \mathbb{P}^2$ in $\mathcal{B}(C)$ with $\delta=\delta(f)$, we assume that $C$ is expressed as 
\[
C=\{ x=[x_1:x_2:x_3] \in \mathbb{P}^2 \mid g(x_1:x_2:x_3)=0 \}, 
\]
where $g$ is a homogeneous polynomial in $\mathbb{Z}[\delta][x_1,x_2,x_3]$, 
and that $f=f_{\delta}$ is also expressed as 
\[
f(x)=f_{\delta}(x)=[f_1(x_1:x_2:x_3):f_2(x_1:x_2:x_3):f_3(x_1:x_2:x_3)] \in \mathbb{P}^2, 
\]
where $f_i$ are homogeneous polynomials in $\mathbb{Z}[\delta][x_1,x_2,x_3]$ with 
$\mathrm{deg}_x f_1=\mathrm{deg}_x f_2=\mathrm{deg}_x f_3$. 
Note that if $\delta \in \mathbb{C}^*$ is an algebraic number, then so is any fixed point $w$ of $f$, 
which enables us to consider Galois conjugates of $\delta$ and $w$, 
and also the eigenvalues $(\mu, \nu)$ of $Df(w)$ are algebraic. 
\begin{proposition} \label{prop:Siegel} 
Under the above assmptions, let $\delta \in \mathbb{C}^*$ be an algebraic number with $|\delta|=1$ that is not a root of unity, and 
$w \in \mathbb{P}^2 \setminus C$ be a fixed point of $f_{\delta}$ outside $C$. 
Moreover assume that there are Galois conjugates $(\delta_*,w_*)$ of $(\delta,w)$ with $|\delta_*|=1$ 
and $f_{\delta_*}(w_*)=w_*$ such that 
\[
\{ \mathrm{Tr}Df_{\delta}(w)\}^2/\mathrm{Det}Df_{\delta}(w) \in [0,4], \qquad
\{ \mathrm{Tr}Df_{\delta_*}(w_*)\}^2/\mathrm{Det}Df_{\delta_*}(w_*) \notin [0,4]. 
\]
Then $f=f_{\delta}$ has a Siegel disk centered at $w$. 
\end{proposition}
{\it Proof}. 
(see \cite{M}). 
Let $(\mu_*, \nu_*)$ be the eigenvalues of $Df_{\delta_*}(w_*)$, which are Galois conjugates of the eigenvalues $(\mu,\nu)$ 
of $Df_{\delta}(w)$. 
Note that $\mu_* \nu_* =\mathrm{Det}  Df_{\delta_*}(w_*)=\delta_*$, as $w_*$ also lies outside $C$.  
Moreover, it should be noted that 
\[
\{ \mathrm{Tr}Df_{\delta}(w)\}^2/\mathrm{Det}Df_{\delta}(w) 
= \frac{(\mu+\nu)^2}{\mu \nu}= \frac{\mu}{\nu} + 
\frac{\nu}{\mu} +2, 
\]
and that a complex number $z \in \mathbb{C}$ satisfies 
$z + z^{-1}+2 \in [0,4]$ if and only if $|z|=1$. 
Hence it follows from our assumption that 
$|\mu/\nu|=1$ and $|\mu_*/\nu_*| \neq 1$. 
Since $|\mu \nu|=|\delta|=1$, we have 
$(\mu, \nu) \in (S^1)^2$. 
Now assume that $\mu^k \nu^{l} =1$ for 
$(k,l) \in \mathbb{Z}^2$. 
Since $(\mu_*, \nu_*)$ are Galois conjugates of $(\mu,\nu)$, one has 
$1=\mu_*^k \nu_*^{l} = (\delta_*)^{(k+l)/2} 
(\mu_*/\nu_*)^{(k-l)/2}$ 
and thus $k= l$ as $|\delta_*|=1$ and $|\mu_*/\nu_*| \neq 1$. 
Since $1 = \mu_*^k \nu_*^k=\delta_*^k$ and 
$\delta_*$ is not a root of unity, we have $k=0$, namely, $(k,l)=(0,0)$. 
Therefore $Df(w)$ is an irrational rotation with the algebraic eigenvalues $(\mu,\nu)$, 
which shows that $f$ has a Siegel disk centered at $w$. 
\hfill $\Box$ \par\medskip 
\section{Birational Maps Preserving a Cuspidal Curve} \label{sec:Quad} 
In this section, we consider a class of quadratic birational maps preserving a cuspidal cubic curve. 
For a parameter $\delta \in \mathbb{C} \setminus \{0,1\}$, 
let us consider a quadratic map $f=f_{\delta} : \mathbb{P}^2 \to \mathbb{P}^2$, 
which is explicitly given by $f [x:y:z]=[f_x:f_y:f_z]$ 
in homogeneous coordinates, where 
\begin{equation} \label{eqn:quad}
\left\{
\begin{array}{ll}
f_x[x:y:z] = & \delta \cdot (x y - 2 d y z +2 d^3 x z -d^4 z^2) \\[2mm]
f_y[x:y:z] = & \delta^3 \cdot (y^2  -3  d^2 x y + 3 d^4 x^2 - d^6 z^2) \\[2mm]
f_z[x:y:z] = & y z - 3 d x^2 + 3 d^2 x z - d^3 z^2 
\end{array}
\right.
\end{equation}
with $d:= (3\delta)^{-1}(1-\delta)$. 
Then $f$ is a birational map preserving the cubic curve 
$C:=\{ yz^2=x^3 \} \subset \mathbb{P}^2$ with a cusp located at $[0:1:0]$, 
and also preserving its smooth locus $C^*=C \setminus \{[0:1:0]\}$. 
Indeed, with the parametrization $p : \mathbb{C} \to C^*$ given by $p(t)=[t:t^3:1]$, 
the restriction of $f$ to $C^*$ is expressed as 
$f|_{C^*} : \mathbb{C} \ni t \mapsto \delta \cdot (t +d) \in \mathbb{C}$. 
The indeterminacy sets of $f^{\pm 1}$ are given by $I(f^{\pm 1})=\{ p_1^{\pm},p_2^{\pm},p_3^{\pm} \}$, 
where $p_1^+:= p(d) \in C^*$ and $p_1^-:= p(- \delta \cdot d) \in C^*$. 
Moreover, for $i=1,2$, the point $p_{i+1}^{\pm}$ is defined by the property 
$\{p_{i+1}^{\pm} \}=C_{i}^{\pm} \cap \mathcal{E}_{i}^{\pm}$, 
where $C_0^{\pm}:=C^*$ and $C_{i}^{\pm}$ is inductively given by the strict transform $(\pi_i^{\pm})^{-1}(C_{i-1}^{\pm})$ 
under the blowup $\pi_i^{\pm}$ of $p_{i}^{\pm}$ with exceptional curve $\mathcal{E}_{i}^{\pm}$. 
In this case, we write $p_1^{\pm} < p_2^{\pm} < p_3^{\pm}$. 
Hence by permitting infinitely near points, we conclude that $I(f^{\pm 1})$ are contained in $C^{*}$, 
and  that $f$ is a quadratic birational map in $\mathcal{B}(C)$ with $\delta(f)=\delta$ 
from the expression for $f|_{C^*}$. 
Conversely, if a quadratic map $f \in \mathcal{B}(C)$ with $I(f)=\{p_1^+,p_2^+,p_3^+\}$ 
satisfies $\delta(f)=\delta$ and 
$p_1^+=p(d)<p_2^+<p_3^+$, then $f=f_{\delta}$ is given by (\ref{eqn:quad}) (see \cite{U1, U2}). 
\par
There are exactly two fixed points $\{w_1, w_2 \}$ of $f$ outside the curve $C$, 
and each point is expressed as $w_{i}=[x_i:r_{\tau}(x_i):1]$, 
where 
\[
r_{\tau}(x):=\frac{\tau-2}{3 (\tau+1)} x - \frac{(\tau-2)^2}{27 (\tau+1)} 
\]
with $\tau:=\delta+1/\delta$, 
and $x_i$ is a root of the quadratic equation 
\[
Q_{\tau}(x) := 27 x^2 - 9 (\tau-2) x +(\tau-1)(\tau-2)=0. 
\]
Moreover  we have 
\[
\frac{\{ \mathrm{Tr}Df(w_i)\}^2}{\mathrm{Det}Df(w_i)} = s(\tau,x_i) := 
\frac{1}{\tau+2} \{ 9(\tau-1) x_i - (\tau^2- 4 \tau +6) \}^2. 
\]
\par
Now in order to construct an automorphism on a rational surface, 
we consider the case where the orbit $p_i^k := f^{k}(p_i^-)$ of each backward indeterminacy point $p_i^-$ reaches 
the forward indeterminacy point $p_i^+$, namely, $p_i^n=p_i^+$ for some $n \ge 1$. 
If such an $n \ge 1$ exists, then 
Proposition \ref{prop:auto} shows that the proper blowup $\pi : X \to \mathbb{P}^2$ for $f$ 
lifts $f $ to an automorphism $F : X \to X$. 
\par
From now on we assume $n=8$. 
As $p_1^k=p(-\delta^{k+1} \cdot d + (1-\delta^k)/3)$, 
it follows from the relation $p(-\delta^{9} \cdot d + (1-\delta^8)/3)=p(d)$ 
that $\delta$ is a root of $(\delta+1)S(\delta)=0$, 
where 
\[
S(\delta)=\delta^8 -2 \delta^7 +\delta^6 -2 \delta^5 + \delta^4 -2 \delta^3 + \delta^2 -2 \delta +1
\]
is a Salem polynomial. 
Conversely, for any root $\delta$ of $S(\delta)=0$, 
the birational map $f=f_{\delta}$ satisfies $p_i^8=p_i^+$ for any $i \in \{1,2,3\}$, 
as $p_1^k<p_2^k<p_3^k$ for any $0 \le k \le 8$, 
and hence lifts to the automorphism $F=F_{\delta} : X \to X$. 
The roots of $S(\delta)=0$ on the real line are $\delta \approx 1.9940, \, 0.5015$ and 
the other roots lie on the unit circle, given by 
$\delta \approx 0.6098 \pm 0.7925 i$, $-0.1098 \pm 0.9939 i$, $-0.7478 \pm 0.6640 i$, which yields 
$\tau \approx 1.2197$, $-0.2197$, $-1.4955$. 
By virtue of Proposition \ref{prop:eigen} (see also the proof of Proposition \ref{pro:sieg}), 
$\lambda \approx 1.9940$ is an eigenvalue of 
$F^* : H^2(X;\mathbb{Z}) \to H^2(X;\mathbb{Z})$ and thus the spectral radius of $F^*$, 
which means that $F$ has positive entropy $h_{\mathrm{top}}(F)= \log \lambda \approx 0.6901 > 0$. 
\par
Now we put $(\delta_0, \tau_0) \approx (0.6098 + 0.7925 i, 1.2197)$ 
and $(\delta_*,\tau_*) \approx (-0.7478 + 0.6640 i, -1.4955)$. 
\begin{lemma} \label{lem:loc}
We have $s(\tau_0,x_i) \in [0,4]$ for any root $x_i$ of $Q_{\tau_0}(x)=0$ and 
$s(\tau_*,x_*) \notin [0,4]$ for some root $x_*$ of $Q_{\tau_*}(x)=0$
\end{lemma}
{\it Proof}. Note that $\tau_0 \in I_0:=[1.219,1.220]$ 
and $\tau_* \in I_*:=[-1.496,-1.495]$. 
Moreover the roots $x_i$ of $Q_{\tau_0}(x)=0$ satisfy either 
$x_i \in I_1 := [0.022, 0.023]$ or $x_i \in I_2 :=[-0.283, -0.282]$ 
as $Q_{\tau}(0.022) < 0$, $Q_{\tau}(0.023) > 0$, 
$Q_{\tau}(-0.283) > 0$, $Q_{\tau}(-0.282) < 0$ for any $\tau \in I_0$, 
and a root $x_*$ of $Q_{\tau_*}(x)=0$ satisfies $x_* \in I_{**}:= [-0.711, -0.710]$ 
as $Q_{\tau}(-0.711) > 0$, $Q_{\tau}(-0.710) < 0$ for any $\tau \in I_*$. 
In particular, we have $s(\tau_0,x_i) \ge 0$ and $s(\tau_*,x_*) \ge 0$. 
A little calculation shows that 
$s(\tau,x) \le s(1.219,0.022) < 2.05 < 4$ 
for any $(\tau,x) \in I_0 \times I_1$, 
$s(\tau,x) \le s(1.220,-0.283) < 3.12 < 4$ 
for any $(\tau,x) \in I_0 \times I_2$ and 
$s(\tau,x) \ge s(-1.495,-0.710) > 5.91 > 4$ 
for any $(\tau,x) \in I_* \times I_{**}$. 
Hence the lemma is established. 
\hfill $\Box$ \par\medskip 
Note that $Q_{\tau_0}(x)$ is irreducible over $\mathbb{Q}[\tau_0]$, and thus 
both $(\delta_0,w_1)$ and $(\delta_0,w_2)$ are Galois conjugates of $(\delta_*,w_*)$. 
Proposition \ref{prop:Siegel} yields the following (see also Remark \ref{rem:iden}). 
\begin{proposition} \label{prop:MI}
The automorphism $F=F_{\delta_{0}} \in \mathcal{QF}(C)$ 
has Siegel disks centered at $w_1, w_2$. 
\end{proposition}
{\it Proof of Theorem \ref{thm:main2}}. 
As $C$ is reduced irreducible, $C$ is either smooth or a nodal cubic or a cuspidal cubic. 
A result of Diller \cite{D} says that there is no automorphism $F \in \mathcal{QF}(C)$ when $C$ is a nodal cubic. 
On the other hand, when $C$ is smooth, the determinant $\delta(F)$ of 
any automorphism $F \in \mathcal{QF}(C)$ is a root of unity. 
Hence for the fixed point $x$, the derivative $DF(x)$ has an eigenvalue $\delta(F)$ if $x \in C$, 
and has the determinant $\mathrm{Det} DF(x)=\delta(F)$ if $x \notin C$. 
In either cases, the eigenvalues of $DF(x)$ are not multiplicatively independent, which means that 
$F$ has no Siegel disk. 
Therefore if $C$ is irreducible and $F \in \mathcal{QF}(C)$ has a Siegel disk, 
then $C$ is a cuspidal cubic curve. 
Moreover, if $C$ is a cuspidal cubic, then $F$ admits at most two Siegel disks by Proposition \ref{pro:sieg}. 
Finally, Proposition \ref{prop:MI} guarantees the existence of the automorphism $F \in \mathcal{QF}(C)$ 
admitting exactly two Siegel disks. 
 \hfill $\Box$ \par
\begin{figure}[t]
\begin{center}
\includegraphics[width=6cm]{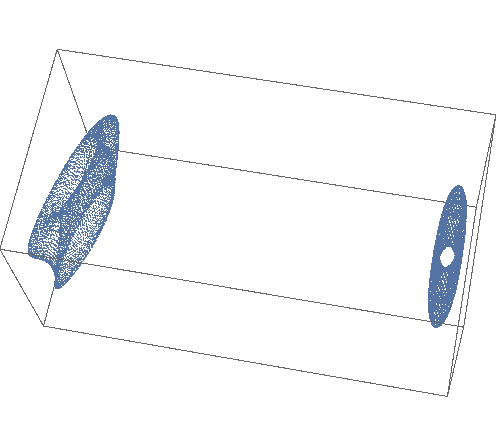}
\end{center}
\caption{Two Siegel disks for automorphism $F$} 
\label{fig:SD01}
\end{figure}

\medskip 
 In Figure \ref{fig:SD01}, we describe two Siegel disks for the automorphism $F$ with the help of Mathematica. 
\section{Birational Maps Preserving Three Lines} \label{sec:lines} 
In this section we consider birational maps preserving three lines meeting a single point. 
To this end, for parameters $\delta \in \mathbb{C}^{\times}$, $a=(a_i)_{i=1}^m \in (\mathbb{C}^{\times})^m$, 
$b=(b_j)_{j=1}^n \in (\mathbb{C}^{\times})^n$, 
let $f =f_{\delta,a,b}: \mathbb{C}^2 \to \mathbb{C}^2$ be a birational map given by
\begin{equation} \label{eqn:birlines}
f : \mathbb{C}^2 \to \mathbb{C}^2, \quad (x,y) \mapsto (f_1(x,y),f_2(x,y))= \Bigl(y, 
\frac{g_1(y)(x+ \delta y)}{\delta\{ (g_2(y)-g_1(y) ) \frac{x}{y} - \delta g_1(y) \}} \Bigr),
\end{equation}
where $g_1(y)=\prod_{i=1}^m (1-y/a_i)$ and $g_2(y) = \prod_{j=1}^n (1-y/b_j)$. 
The map $f$ preserves the three lines $C= L_1 \cup L_2 \cup L_3$, 
where $L_1=\{x=0\}$, $L_2=\{x+ \delta y=0\}$, $L_3=\{ y=0 \}$, 
and sends these lines as 
\begin{equation} \label{eqn:lines}
f|_{L_1}(0,y) = (y, - \frac{y}{\delta}) \in L_2, \quad f|_{L_2} (- \delta y,y) = (y,0) \in L_3, \quad
f|_{L_3} (x,0) = (0, \frac{- x}{\delta^2 + \delta c x}) \in L_1. 
\end{equation}
Here and hereafter, we use the following notations: 
\begin{equation} \label{eqn:par}
\begin{array}{ll}
\displaystyle \alpha = \sum_{i=1}^m \frac{1}{a_i}, & \displaystyle \beta = \sum_{j=1}^n \frac{1}{b_j}, \\[2mm]
\displaystyle \alpha_0:=\prod_{i=1}^m \frac{1}{a_i}, & \displaystyle \beta_0:=\prod_{j=1}^n \frac{1}{b_j}, \\[2mm]
c= \beta - \alpha. & ~
\end{array}
\end{equation}
Note that the map (\ref{eqn:birlines}) is derived under a certain assumption as in the following lemma. 
\begin{lemma} \label{lem:thli}
Assume that a birational map $h : \mathbb{C}^2 \to \mathbb{C}^2$ of the form $h(x,y)=(y,h_2(x,y))$ satisfies 
$h(L_i)=L_{i+1}$ for $i=1,2,3~(\mathrm{mod}~3)$. Then we have $h=f_{\delta,a,b}$ for some $\delta$, $a=(a_i)$ and $b=(b_j)$.  
\end{lemma}
{\it Proof}. Since $h$ is a birational map, for a generic $(x_0,y_0) \in \mathbb{C}^2$, 
the equation $h(x,y)=(y,h_2(x,y))=(x_0,y_0)$, or $h_2(x,x_0)=y_0$ has a unique root for $x$. 
Hence $h_2(x,y)$ is a rational function of degree $1$ with respect to $x$. 
As $h_2(-\delta y,y)=0$, $h_2(0,y)=-y/\delta$ and $h_2(x,0) \neq 0$, $h_2$ has the form 
$h_2(x,y)=g_1(y) (x+ \delta y) / (g_3(y) x - \delta^2 g_1(y))$ with $g_1(0) \neq 0$. 
By multiplying the denominator and numerator by a common constant if necessary,  one can put 
$g_1(y)= \prod_{i=1}^m(1-y/a_i)$ and then $g_2(y)= g_1(y)+ y g_3(y) / \delta  = \prod_{j=1}^n(1-y/b_j)$, 
which yields the lemma. 
\hfill $\Box$ \par\medskip 
From now on, we assume the following: 
\begin{assumption} \label{ass1} $m=n=N$. 
\end{assumption}
With the embedding $\mathbb{C}^2 \ni (x,y) \hookrightarrow [x:y:1] \in \mathbb{P}^2$, 
we will regard the birational map $f$ and the lines $C$ as those on $\mathbb{P}^2$. 
Then the indeterminacy sets of $f^{\pm 1}$ are given by 
$I(f^{\pm 1})=\{p_{a,i}^{\pm} \}_{i=1}^N \cup \{ p_{b,j}^{\pm} \}_{j=1}^N \cup \{p_{0}^{\pm} \}$,
where 
\begin{equation*}
\begin{array}{lll}
p_{a,i}^{+} = [0:a_i:1], & p_{b,j}^{+}=[-b_j \delta : b_j : 1], & p_0^+=[1:0:0], \\[2mm]
p_{a,i}^{-} = [a_i:0:1], & p_{b,j}^{-}=[b_j : - b_j/\delta : 1], & p_0^-=[0:1:0]. 
\end{array} 
\end{equation*}
Since any indeterminacy point of $f^{\pm 1}$ lies on the smooth locus $C^*$ of the three lines $C$, 
we can conclude that $f \in \mathcal{B}(C)$. 
Moreover, it follows from (\ref{eqn:lines}) that $\delta=\delta(f)$ is the determinant of $f$.  
\begin{remark} \label{rem:curve}
The birational map $f$ contracts curves to indeterminacy points as follows: 
\[
\begin{array}{l}
L_{i}^a:=\{ [x:a_i:1] \, | \, x \in \mathbb{P}^1 \} \to p_{a,i}^-, \\[2mm]
L_{j}^b:=\{ [x:b_j:1] \, | \, x \in \mathbb{P}^1 \} \to p_{b,j}^-, \\[2mm]
D:=\{ [x:y:1] \, | \, (g_2(y)-g_1(y) )x/y -\delta g_1(y)=0 \} \to p_0^-. 
\end{array}
\]
The curves $L_i^a$ and $L_j^b$ are lines passing through $\{p_{a,i}^+, p_0^+ \}$ and 
$\{p_{b,j}^+, p_0^+ \}$ respectively, 
and $D$ is a curve of degree $N$ passing through 
$I(f)$ with multiplicities $\mathrm{mult}_{p_{a,i}^+} D= \mathrm{mult}_{p_{b,j}^+} D=1$ and 
$\mathrm{mult}_{p_0^+} D=N-1$. 
A straightforward calculation shows that the blowup of $p_0^-$ lifts $f$ to a birational map whose restriction to 
$D$ is an isomorphism to the exceptional curve of the blowup. 
Similarly, if $a_i \neq a_k$ for any $k \neq i$, then the blowup of $p_{a,i}^-$ lifts $f$ to a birational map 
whose restriction to $L_i^a$ is an isomorphism to the exceptional curve, 
and also if $b_j \neq b_k$ for any $k \neq j$, then the blowup of $p_{b,j}^-$ lifts $f$ to a birational map 
whose restriction to $L_j^b$ is an isomorphism to the exceptional curve. 
Moreover, the pullback of a generic line by $f$ is a curve $\mathcal{D}$ of degree $N+1$ passing through 
$I(f)$ with multiplicities $\mathrm{mult}_{p_{a,i}^+} \mathcal{D}= \mathrm{mult}_{p_{b,j}^+} \mathcal{D}=1$ and 
$\mathrm{mult}_{p_0^+} \mathcal{D}=N$. 
\end{remark}
Next we determine the fixed points of $f : \mathbb{P}^2 \to \mathbb{P}^2$.  
The fixed points of $f$ on $\mathbb{C}^2$ are given by the singular point $(0,0)$ of $C$, and 
$(x_l,x_l) \in \mathbb{C}^2$, where $x_l$ are the roots of the equation
\begin{equation} \label{eqn:fx1}
\frac{(1+\delta)^2}{\delta} \prod_{i=1}^N (1-\frac{x_l}{a_i}) = \prod_{j=1}^N (1-\frac{x_l}{b_j}). 
\end{equation}
Moreover under Assumption \ref{ass1}, the birational map $f : \mathbb{P}^2 \to \mathbb{P}^2$ preserves the line 
$L=\{[x:y:z] \, | \, z=0\}$ at infinity, and the restriction $f|_L$ is expressed as 
\begin{equation} \label{eqn:fai}
f[x:y:0]=[\delta (\beta_0-\alpha_0)x-\delta^2 \alpha_0 y  : \alpha_0(x+\delta y) : 0],
\end{equation}
where $\alpha_0$, $\beta_0$ are given in (\ref{eqn:par}). 
Hence the fixed points of $f : \mathbb{P}^2 \to \mathbb{P}^2$ lying on $L$ are given by $[x_l:1:0]$, 
where $x_l$ are the roots of the equation 
\begin{equation} \label{eqn:fx2}
\alpha_0 x_l^2 + \delta (2 \alpha_0-\beta_0) x_l + \alpha_0 \delta^2=0.
\end{equation}
Consequently, we have
\begin{proposition} \label{prop:allfix}
The fixed points of $f : \mathbb{P}^2 \to \mathbb{P}^2$ are given by $w_0=[0:0:1] \in C$, 
$w_l=[x_l:x_l:1] \in \mathbb{C}^2$ for $l \in \{ 1,\dots,N \}$, where $x_l$ are the roots of (\ref{eqn:fx1}), 
and $w_l=[x_l:1:0] \in L$ for $l \in \{N+1, N+2 \}$, where $x_l$ are the roots of (\ref{eqn:fx2}). 
Moreover, when $l \in \{ 1, \dots, N+2\}$, the fixed point 
$w_l$ lies outside $C$ and hence satisfies $\mathrm{Det} Df(w_l)=\delta$. 
\end{proposition}
\begin{remark} \label{rem:sing}
It is straightforward to calculate that 
the eigenvalues of $Df(w_0)$ at the singular point $w_0$ of $C$ 
are given by $(\omega \delta^{-1},\omega^{-1} \delta^{-1})$, 
where $\omega$ is a primitive cube root of unity. 
Therefore a Siegel disk is never centered at $w_0$, 
as $(\omega \delta^{-1},\omega^{-1} \delta^{-1})$ are not multiplicatively independent. 
\end{remark}
Now, for $\mathbb{R}^*:=\mathbb{R} \setminus \{ 0 \}$, we put 
\[
A:=\{ c=(\delta,a,b) \in S^1 \times (\mathbb{R}^*)^N \times (\mathbb{R}^*)^N \, | \, 
a_i \neq a_j, b_i \neq b_j, (i \neq j),~ \sum_{j=1}^N \frac{1}{b_j} - \sum_{i=1}^N \frac{1}{a_i} = 1 \}, 
\]
and for $c_0=(\delta_0,a_0,b_0) \in A$ and $\varepsilon >0$, put  
\[
A(c_0 ; \varepsilon):=\{ (\delta,a,b) \in A \, | \, 
|\delta-\delta_0|< \varepsilon, |a - a_0|< \varepsilon, |b-b_0|< \varepsilon \}. 
\]
Then we have the following proposition, whose proof is given in Section \ref{sec:pp2}. 
\begin{proposition} \label{prop:fixrel}
Under the above notations, there exists $\varepsilon>0$ and $c_0, c_* \in A$ such that 
\begin{enumerate}
\item $\displaystyle \frac{\{\mathrm{Tr} Df(w_l)\}^2}{\mathrm{Det} Df(w_l)} \in [0,4]$ for any 
$l \in \{1, \dots ,N+2\}$ if $(\delta,a,b) \in A(c_0;\varepsilon)$, 
\item $\displaystyle \frac{\{\mathrm{Tr} Df(w_l)\}^2}{\mathrm{Det} Df(w_l)} \notin [0,4]$ for any 
$l \in \{1, \dots ,N+2\}$ if $(\delta,a,b) \in A(c_*;\varepsilon)$.
\end{enumerate}
\end{proposition}
\par
It should be noted that the indeterminacy point $p_0^- \in I(f^{-1})$ satisfies 
$f^2(p_0^-)=p_0^+ \in I(f)$. 
Furthermore we assume the following: 
\begin{assumption} \label{ass2}
For given parameters $m=(m_i)_{i=1}^N, n=(n_j)_{j=1}^N \in \mathbb{N}^N$ 
except for $(m,n)=((1),(1)) \in (\mathbb{N}^1)^2$, the map $f : \mathbb{P}^2 \to \mathbb{P}^2$ satisfies  
\begin{equation} \label{eqn:ab}
\begin{array}{rl}
f^{3m_i-2}(p_{a,i}^-) & = p_{a,i}^+ \quad (i=1,\dots,N), \\ 
f^{3 n_j}(p_{b,j}^-) & = p_{b,j}^+  \quad (j=1,\dots,N).
\end{array}
\end{equation}
\end{assumption}
\begin{lemma} \label{lem:ab}
Under Assumption 2, we have 
\[
\frac{1}{a_i} = - \frac{\delta (\delta^{3 m_i}-1)}{(\delta^3-1) (\delta^{3m_i-1} +1)} c, \qquad 
\frac{1}{b_j} = \frac{\delta^2 (\delta^{3 n_j}-1)}{(\delta^3-1) (\delta^{3n_j+1} +1)} c, 
\]
where $c=\beta-\alpha$ is given in (\ref{eqn:par}). 
In particular, if $c \neq 0$, then $\delta$ satisfies the equation 
\begin{equation} \label{eqn:delta}
\chi_{m,n}(\delta) := \sum_{j=1}^N \frac{\delta^2 (\delta^{3 n_j}-1)}{(\delta^3-1) (\delta^{3n_j+1} +1)} 
+ \sum_{i=1}^N \frac{\delta (\delta^{3 m_i}-1)}{(\delta^3-1) (\delta^{3m_i-1} +1)} = 1. 
\end{equation}
\end{lemma}
{\it Proof}. It follows from (\ref{eqn:lines}) that $f^3(0,y)=(0,h_1(y))$, $f^3(x,-x/\delta)=(h_1(x),-h_1(x)/\delta)$ and 
hence $f^{3k}(0,y)=(0,h_k(y))$, $f^{3k}(x,-x/\delta)=(h_k(x),-h_k(x)/\delta)$, where 
\[
h_k(x):=\frac{1}{\delta^{3 k}(\frac{1}{x} -p)+p}, \qquad 
p:=\frac{\delta c}{(\delta^3-1)}. 
\]
Since $f(a_i,0)=(0, -a_i \{ \delta (\delta+  c a_i)\}^{-1} )$, the assumption (\ref{eqn:ab}) is equivalent to 
$h_{m_i-1}(-a_i \{ \delta (\delta+  c a_i)\}^{-1})=a_i$ and $h_{n_j}(b_j)=-b_j \delta$, 
which yield the desired expressions for $1/a_i$ and $1/b_j$. 
Finally, the relation (\ref{eqn:delta}) follows from $c=\beta-\alpha=\sum_{j=1}^N1/b_j - \sum_{i=1}^N1/a_i$. 
\hfill $\Box$ \par\medskip 
Conversely, for given $m=(m_i), n=(n_j) \in \mathbb{N}^N$, let $\delta \in \mathbb{C}^*$ be any root of (\ref{eqn:delta}), 
and $a=(a_i)$, $b=(b_j)$ be parameters given by $a_i=a_{m_i}(\delta)$, $b_j=b_{n_j}(\delta)$, where 
\begin{equation} \label{eqn:abk}
a_k(\delta) := - \frac{(\delta^3-1) (\delta^{3k-1}+1)}{\delta (\delta^{3k}-1)}, \qquad 
b_k(\delta) := \frac{(\delta^3-1) (\delta^{3k+1}+1)}{\delta^2 (\delta^{3k}-1)}. 
\end{equation}
Then the birational map $f=f_{\delta,a,b}$ satisfies the condition (\ref{eqn:ab}). 
Proposition \ref{prop:auto} shows that there is a proper blowup $\pi : X \to \mathbb{P}^2$ for $f$, 
and $\pi$ lifts $f : \mathbb{P}^2 \to \mathbb{P}^2$ to an automorphism $F_{m,n} : X \to X$. 
Note that the points blown up by $\pi$ satisfy the assumptions in Proposition \ref{prop:eigen}.
Thus the root $\delta$ of the equation (\ref{eqn:delta}), which is the determinant of $f$, is an eigenvalue of 
$F_{m,n}^* : H^2(X;\mathbb{Z}) \to H^2(X;\mathbb{Z})$. 
On the other hand, under Assumption \ref{ass2}, 
there exists $\lambda >1$ so that $\chi_{m,n}(\lambda)=1$
since $\chi_{m,n}(1)>1$ and $\lim_{\delta \to \infty} \chi_{m,n}(\delta)=0$. 
Hence $\lambda=\lambda_{m,n}:=\lambda(F_{m,n}^*)>1$ is the spectral radius, 
which is a root of the Salem polynomial $S_{m,n}(t):=S_{F_{m,n}}(t)=0$. 
As $S_{m,n}(t)$ is irreducible, any root of $S_{m,n}(t)=0$ is a root of $\chi_{m,n}(t)=1$. 
Therefore we have
\begin{corollary} \label{cor:salem} 
Under the assumption that $(m,n) \neq ((1),(1))$, any root $\delta$ of $S_{m,n}(t)=0$ satisfies $\chi_{m,n}(\delta)=1$. 
Moreover, the birational map $f=f_{\delta,(a_{m_i}(\delta)),(b_{n_j}(\delta))}$ lifts to the automorphism $F_{m,n}$, 
having positive entropy 
$h_{\mathrm{top}}(F_{m,n})= \mathrm{log} \lambda_{m,n}>0$ with the spectral radius $\lambda_{m,n}=\lambda(F_{m,n}^*)>1$. 
\end{corollary} 
\begin{lemma} \label{lem:dense}
If $\delta \in S^1$ is given by $\delta=\exp(2 \pi i \nu)$ with an irrational real number $\nu$, then 
$\{a_{k}(\delta) \}_{k \in \mathbb{N}}$ and $\{b_{k}(\delta) \}_{k \in \mathbb{N}}$ 
are sequences of real numbers and dense in $\mathbb{R}$. 
\end{lemma}
{\it Proof}. First we notice that 
\[\begin{array}{rl}
\displaystyle -\frac{(\delta^3-1) (\delta^{3k-1}+1)}{\delta (\delta^{3k}-1)} = & \displaystyle 
- \frac{(\delta^{3/2}-\delta^{-3/2}) (\delta^{(3k-1)/2}+\delta^{-(3k-1)/2})}{(\delta^{3k/2}-\delta^{-3k/2})} \\ \displaystyle 
= &\displaystyle - 2 \frac{\sin(3 \pi \nu) \cos \{ (3 k-1) \pi \nu\} }{\sin (3 k \pi \nu)} 
= -2 \sin(3 \pi \nu) \Bigl\{\frac{\cos (\pi \nu)}{\tan(3k \pi \nu)} + \sin(\pi \nu) \Bigr\}, \\
\displaystyle \frac{(\delta^3-1) (\delta^{3k+1}+1)}{\delta^2 (\delta^{3k}-1)} = & \displaystyle 
\frac{(\delta^{3/2}-\delta^{-3/2}) (\delta^{(3k+1)/2}+\delta^{-(3k+1)/2})}{(\delta^{3k/2}-\delta^{-3k/2})} \\ \displaystyle 
= &\displaystyle 2 \frac{\sin(3 \pi \nu) \cos \{ (3 k+1) \pi \nu\} }{\sin (3 k \pi \nu)} 
= 2 \sin(3 \pi \nu) \Bigl\{\frac{\cos (\pi \nu)}{\tan(3k \pi \nu)} - \sin(\pi \nu) \Bigr\}
\end{array}
\]
are real numbers. 
Since $\{3 k \pi \nu \}_{k=1}^\infty \subset (- \pi/2,\pi/2)~(\mathrm{mod}~\pi)$ is dense from the irrationality of $\nu$, 
so is $\{ \tan(3 k \pi \nu) \}_{k=1}^\infty \subset \mathbb{R}$, which establishes the lemma 
as $\sin(3 \pi \nu) \cos(\pi \nu) \neq 0$.  
\hfill $\Box$ \par\medskip 
\begin{proposition} \label{prop:equi}
The roots of $S_{m,n}(t)=0$ other than $\lambda_{m,n}^{\pm 1}$ are equidistributed on $S^1$ 
as either $m_i \to \infty$ for some $i$ or $n_j \to \infty$ for some $j$.  
\end{proposition}
{\it Proof}. 
A result of Bilu (see \cite{B, M}) says that if $\{ \rho_k\}_{k \in \mathbb{N}}$ is a sequence of algebraic units with 
$\lim_{k \to \infty} \mathrm{deg}(\rho_k) =\infty$ then $\{ \overline{\delta}_{\rho_k} \}$ weakly 
converges to the normalized Haar measure on $S^1$. 
Here for an algebraic number $\rho \neq 0$, we put 
\[
\overline{\delta}_{\rho} := \frac{1}{\mathrm{deg}(\rho)} \sum_{\rho' \underset{\mathrm{conj.}}{\sim} \rho} \delta_{\rho'}
\]
with the Dirac measure $\delta_{\rho'}$ at $\rho'$. 
Since $\lambda_{m,n}$ satisfies $\lambda_{m,n} \to \lambda< \infty$ as $m_i \to \infty$ or $n_j \to \infty$ and $\lambda$ is not 
a Salem number, we have $\mathrm{deg}(\lambda_{m,n}) \to \infty$. 
As $\lambda_{m,n}$ is an algebraic unit, the proposition is established. 
\hfill $\Box$ \par\medskip 
\begin{proposition} \label{prop:approx}
Let $\varepsilon>0$ and $c_0, c_* \in A$ be given in Proposition \ref{prop:fixrel}, 
and  $a_k(\delta)$, $b_k(\delta)$ be given in (\ref{eqn:abk}). 
Then there exist $m,n \in \mathbb{N}^N$ and $\delta_0, \delta_* \in S^1$ such that 
\begin{enumerate}
\item $S_{m,n}(\delta_0)=S_{m,n}(\delta_*)=0$, and
\item $(\delta_0,(a_{m_i}(\delta_0)),(b_{n_j}(\delta_0))) \in A(c_0;\varepsilon)$, \quad
$(\delta_*,(a_{m_i}(\delta_*)),(b_{n_j}(\delta_*))) \in A(c_*;\varepsilon)$. 
\end{enumerate}
\end{proposition}
{\it Proof}. We put $c_0=(d_0,(a_i^0),(b_j^0))$, $c_*=(d_*,(a_i^*),(b_j^*))$, and 
without loss of generality, we may assume that $d_0$ and $d_*$ are multiplicatively independent. 
Then from Lemma \ref{lem:dense}, one can fix $(m_i)_{i=1}^{N-1}$ and $(n_j)_{j=1}^N$ so that 
$a_i^0 \approx a_{m_i}(d_0)$, $a_i^* \approx a_{m_i}(d_*)$ for $i \in \{1,\dots,N-1 \}$ and 
$b_j^0 \approx b_{n_j}(d_0)$, $b_j^* \approx b_{n_j}(d_*)$ for $j \in \{1,\dots,N\}$. 
By Proposition \ref{prop:equi}, there exists $m_N >>1$ such that roots $\delta_0, \delta_* \in S^1$ of 
$S_{m,n}(t)=0$ satisfy $\delta_0 \approx d_0$, $\delta_* \approx d_*$ and hence 
$a_i^0 \approx a_{m_i}(\delta_0)$, $a_i^* \approx a_{m_i}(\delta_*)$ for $i \in \{1,\dots,N-1 \}$ and 
$b_j^0 \approx b_{n_j}(\delta_0)$, $b_j^* \approx b_{n_j}(\delta_*)$ for $j \in \{1,\dots,N\}$. 
As 
\[
\sum_{j=1}^N \frac{1}{b_j} - \sum_{i=1}^N \frac{1}{a_i} = \chi_{m,n}(\delta_0)=\chi_{m,n}(\delta_*)=1
\]
from Corollary \ref{cor:salem}, 
we have $a_N^0 \approx a_{m_N}(\delta_0)$, $a_N^* \approx a_{m_N}(\delta_*)$ so that 
the condition (2) holds. 
\hfill $\Box$ \par\medskip 
For the parameters given in Proposition \ref{prop:approx}, 
fix the birational maps $f_0=f_{\delta_0,(a_{m_i}(\delta_0)),(b_{n_j}(\delta_0))}$ and 
$f_*=f_{\delta_*,(a_{m_i}(\delta_*)),(b_{n_j}(\delta_*))}$. 
As $f_0$ and $f_*$ are Galois conjugate and each fixed point of $f_0$ outside $C$ is a Galois conjugate of 
a fixed point of $f_*$ outside $C$, Propositions \ref{prop:Siegel},  \ref{prop:fixrel}, \ref{prop:approx}
yield the following corollary. 
\begin{corollary} \label{cor:siegel} 
The birational map $f_0$ has $N+2$ fixed points $w_1,\dots,w_{N+2}$ at which Siegel disks are centered. 
\end{corollary} 
\begin{proposition} \label{pro:nth}
Let $F : X \to X$ be the automorphism that is the lift of $f_0$ by the proper blowup $\pi : X \to \mathbb{P}^2$ for $f_0$. 
Then, $F$ has positive entropy $h_{\mathrm{top}}(F)= \mathrm{log} \lambda_{m,n}>0$ and 
has exactly $N+3$ isolated fixed points $w_0,\dots,w_{N+2}$ (see also Remark \ref{rem:iden}). 
\end{proposition}
{\it Proof}. 
Corollary \ref{cor:salem} says that $F=F_{m,n}$ has positive entropy $h_{\mathrm{top}}(F_{m,n})= \mathrm{log} \lambda_{m,n}>0$. 
Now note that the indeterminacy points $I(f^{\pm 1})$ are blown up by $\pi$. 
Remark \ref{rem:curve} says that $F^*$ sends curves as 
\[
\begin{array}{l}
[H] \longmapsto (N+1) [H] - N[E_0^+] - \sum_{i=1}^N [E_{a,i}^+] - \sum_{j=1}^N [E_{b,j}^+], \\[2mm]
[E_0^-] \longmapsto N [H] - (N-1) [E_0^+] - \sum_{i=1}^N [E_{a,i}^+] - \sum_{j=1}^N [E_{b,j}^+], \\[2mm]
[E_{a,i}^-] \longmapsto [H] - [E_0^+] - [E_{a,i}^+], \\[2mm]
[E_{b,j}^-] \longmapsto  [H] - [E_0^+] - [E_{b,j}^+], 
\end{array}
\]
where $E_0^{\pm}, E_{a,i}^{\pm}, E_{b,j}^{\pm}$ are the exceptional divisors over the points 
$p_0^{\pm}, p_{a,i}^{\pm}, p_{b,j}^{\pm}$, respectively. 
It follows from Lemma \ref{lem:exchange} that any exceptional divisor over the point outside $I(f^{-1})$ is sent 
to another exceptional one by $F^*$. 
Hence we have $\mathrm{Tr}(F^*|_{H^2(X;\mathbb{Z})}) \le N+1$. 
Proposition \ref{pro:nof} says that there are at most $N+3$ isolated fixed points for $F$, 
and the existence of the fixed points $w_0, \dots, w_{N+2}$ given in Proposition \ref{prop:allfix} says that there are exactly 
$N+3$ isolated fixed points for $F$. 
\hfill $\Box$ \\[2mm] 
{\it Proof of Theorem \ref{thm:main1}}. 
First assume $k \ge 3$ and put $N=k-2$. 
The automorphism $F$ mentioned in Proposition \ref{pro:nth} has positive entropy and 
has exactly $k+1$ fixed points $w_0,\dots,w_k$. 
Among the fixed points, no Siegel disk is centered at $w_0$ from Remark \ref{rem:sing}, 
and Siegel disks are centered at $w_1,\dots,w_k$ from Corollary \ref{cor:siegel}. 
Therefore $F$ is a desired automorphism satisfying the condition mentioned in Theorem \ref{thm:main1}. 
\par
When $k=0,1$, McMullen \cite{M} and Bedford-Kim \cite{BK1} showed the existence of an automorphism $F$ 
satisfying the condition. 
The automorphism $F$ realizes the Coxeter element and is obtained from a birational map 
$f : \mathbb{P}^2 \to \mathbb{P}^2$ of degree $2$ by blowing up points on the smooth locus of a cubic curve $C$. 
Moreover, $C$ is a cuspidal cubic if $k=0$, and $C$ is either conic with a tangent line or three lines through a point 
if $k=1$. 
Finally, when $k=2$, the existence is shown in Theorem \ref{thm:main2}. 
The theorem is established. 
\hfill $\Box$ \par
\begin{figure}[t]
\begin{center}
\includegraphics[width=8cm]{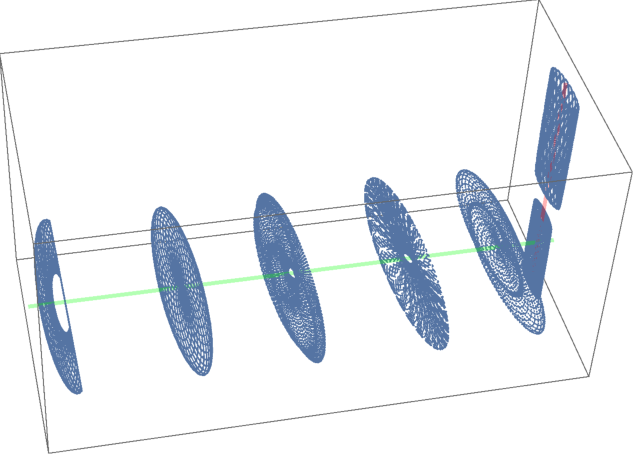}
\end{center}
\caption{Siegel disks for an automorphism ($N=5$)} 
\label{fig:SD02}
\end{figure}

\medskip 
 With the help of Mathematica, we describe Siegel disks of an automorphism for the parameters $N=5$, $m=(280,104,54,36,27)$, $n=(205,381,432,450,459)$ and $\delta \approx -0.5037+ 0.8639 i $ in Figure \ref{fig:SD02}. 
\section{Proof of Proposition \ref{prop:fpf}} \label{sec:pp1} 
This section is devoted to the proof of Proposition \ref{prop:fpf}. 
Since automorphisms may fix a curve pointwise, 
we use S. Saito's fixed point formula instead of a classical fixed point one (see \cite{IU, S}). 
Let $X$ be a smooth projective surface and $f : X \to X$ be an automorphism different from the identity. 
Then the idea of his fixed point formula is to divide the set $X_1(f)$ of irreducible curves fixed pointwise by $f$ 
into the curves of type I and those of type I\!I: 
\[
X_1(f) = X_I(f) \amalg X_{I\!I}(f), 
\]
and to contribute different types of curves to the formula in different ways. 
Namely, the formula says that the Lefschetz number 
\[L(f) := \sum_{i} (-1)^i \, \mathrm{Tr} 
[\, f^{*} : H^i(X;\mathbb{Z}) \to H^i(X;\mathbb{Z})\,] 
\]
of the automorphism $f$ is expressed as 
\[
L(f) = \sum_{x \in X_0(f)} \nu_x(f) + 
\sum_{C \in X_{I}(f)} \chi_C \cdot \nu_C(f) + 
\sum_{C \in X_{I\!I}(f)} \tau_C \cdot \nu_C(f), 
\]
where $X_0(f)$ is the set of fixed points of $f$, $\chi_C$ is the Euler characteristic of the normalization 
of $C \in X_I(f)$ and $\tau_C$ is the self-intersection number 
of $C \in X_{I\!I}(f)$. 
We shall omit the precise definitions of the indices $\nu_x(f)$ and $\nu_{C}(f)$. 
However, it is known that $\nu_{C}(f)$ is a positive integer, 
and $\nu_x(f)$ is a nonnegative integer, which is positive if $x \in X_0(f)$ is an isolated fixed point. 
On the other hand, the types of fixed curves are defined by using the action of $f$ on the completion $A_x$ of 
the local ring of $X$ at $x$, which is isomorphic to the formal power series ring $\mathbb{C}[\![z_1,z_2]\!]$, 
as $X$ is assumed to be a smooth surface. 
Now given a fixed curve $C \in X_1(f)$, we take a smooth 
point $x$ of $C$ and identify $A_x$ with $\mathbb{C}[\![z_1,z_2]\!]$ in 
such a manner that $C$ has the local defining equation $z_1 = 0$ near $x$. 
Then the induced automorphism $f_x^* : A_x \to A_x$ can be expressed as 
\begin{equation} \label{eqn:act}
\left\{
\begin{array}{rcl}
f_x^*(z_1) &=& z_1 + z_1^k \cdot f_1 \\[2mm]
f_x^*(z_2) &=& z_2 + z_1^l \cdot f_2
\end{array}
\right. 
\end{equation}
for some $k$, $l \in \mathbb{N} \cup \{\infty\}$ and some $f_i \in A_x$ 
such that $f_i(0,z_2)$ is a nonzero element of $\mathbb{C}[\![z_2]\!]$. 
Here we put $z_1^{\infty} := 0$ by convention. 
Then it turns out (see \cite{IU}, Lemma 6.1) that $\nu_{C}(f) = \min \{ k,l \}$ and $C \in X_{I}(f)$ if and only if $k \le l$, 
which is independent of the choice of the smooth point $x$ on $C$ and the coordinates $z_1, z_2$. 
Note that if the derivative $Df(x)$ has an eigenvalue different from $1$, then the relation (\ref{eqn:act}) 
yields $k=1$ and $f_1(0,0) \neq 0$. 
In particular, the fixed curve $C$ must be of type $I$. 
\\[2mm]
{\it Proof of Proposition \ref{prop:fpf}}. 
Now if $X$ is a rational surface, then the cohomology group of $X$ is expressed as 
\[
H^i(X;\mathbb{Z}) \cong 
\left\{
\begin{array}{ll}
\mathbb{Z}^{\rho+1} & (i=2) \\
\mathbb{Z} & (i=0,4) \\
0 & (i \neq 0,2,4) \\
\end{array}
\right. 
\]
for some $\rho \ge 0$. 
Moreover, if $F$ is an automorphism on $X$, then the action $F^*$ on $H^i(X;\mathbb{Z})$ is trivial for $i=0,4$, 
which shows that $L(F)=\mathrm{Tr} (F^*|_{H^2(X;\mathbb{Z})})+2$. 
On the other hand, the above argument says that any fixed curve is of type $I$. 
Furthermore if $F$ has positive entropy, then it is known (see \cite{DJS}) 
that any fixed curve $C$ has nonnegative Euler characteristic $\chi_C \ge 0$. 
Hence the fixed point formula says that $F$ has at most $L(F)=\mathrm{Tr} (F^*|_{H^2(X;\mathbb{Z})})+2$ 
isolated fixed points.  
\hfill $\Box$ \par\medskip 
\section{Proof of Proposition \ref{prop:fixrel}} \label{sec:pp2} 
In this section, we will prove Proposition \ref{prop:fixrel}. 
To this end we need some auxiliary lemmas. 
Let $f$ be the birational map given by (\ref{eqn:birlines}) with $m=n=N$. 
\begin{lemma} \label{lem:fxtr1}
For any fixed point $(x_l,x_l) \in \mathbb{C}^2$ with $x_{l}$ satisfying (\ref{eqn:fx1}), we have 
\[
\mathrm{Tr} Df (x_l,x_l)=\frac{\partial f_2}{\partial y}(x_l,x_l)= (\delta+1) \Bigl\{ 1
- \sum_{i=1}^N \frac{1}{1-x_l/a_i} + \sum_{j=1}^N \frac{1}{1-x_l/b_j} \Bigr\}. 
\]
\end{lemma}
{\it Proof}. First it follows from $(f_1)_x=0$ that 
$\mathrm{Tr} Df (x_l,x_l)=(f_2)_y(x_l,x_l)$. 
Moreover, by the relation $g_2(x_l)=g_1(x_l)(1+\delta)^2 /\delta$, one has 
\[
\frac{\partial f_2}{\partial y}(x_l,x_l)= 
\frac{~g_1(x_l) (1+\delta)^2 + x_l g_1'(x_l) (1+\delta)^2 - x_l g_2'(x_l)\delta~}{g_1(x_l) (1+\delta)}. 
\]
Therefore by combining the relations
\[
\begin{array}{l}
\displaystyle 
x_l g_1'(x_l)= g_1(x_l) \sum_{i=1}^N \frac{-x_l/a_i}{1-x_l/a_i}= g_1(x_l) \{N - \sum_{i=1}^N \frac{1}{1-x_l/a_i} \}, \\[2mm]
\displaystyle
x_l g_2'(x_l) \delta= g_2(x_l) \delta \sum_{j=1}^N \frac{-x_l/b_j}{1-x_l/b_j}=g_1(x_l) (1+\delta)^2 \{N - \sum_{j=1}^N \frac{1}{1-x_l/b_j} \}, 
\end{array}
\]
we obtain the desired from. 
\hfill $\Box$ \par\medskip 
\begin{lemma} \label{lem:fxtr2}
Assume $\delta \in S^1$. 
For any fixed point $w_l=[x_l:1:0] \in L$ with $x_{l}$ satisfying (\ref{eqn:fx2}), we have 
\[
\frac{\{\mathrm{Tr} Df(w_l)\}^2}{\mathrm{Det} Df(w_l)}  \in [0,4]  \Longleftrightarrow 
\frac{\beta_0}{\alpha_0} \in [0,4]. 
\]
\end{lemma}
{\it Proof}. We use the fact that the eigenvalues of $Df$ at $w_l=[x_l:1:0]$ for $l \in \{N+1,N+2\}$ 
are given by $(\delta x_l^{-1},x_l)$. 
It follows from the equation (\ref{eqn:fx2}) that $t:=\delta^{-1} x_l$ satisfies 
\[
t= \frac{1}{2} \Bigr\{ \frac{\beta_0}{\alpha_0} - 2 \pm 
\sqrt{\frac{\beta_0}{\alpha_0} \Bigl(\frac{\beta_0}{\alpha_0} - 4 \Bigr)} \Bigr\}. 
\]
Moreover, one has $\{ \mathrm{Tr}Df(w_l)\}^2/\mathrm{Det}Df(w_l)= 2+ \delta t^2 + (\delta t^2)^{-1}$.  
As $\delta \in S^1$, it turns out that $\{ \mathrm{Tr}Df(w_l)\}^2/\mathrm{Det}Df(w_l) \in [0,4]$ if and only if 
$t \in S^1$, or in other words, $\beta_0/\alpha_0 \in [0,4]$. 
\hfill $\Box$ \par\medskip 
Now we show the existence of the parameters $c,c_* \in A$ mentioned in Proposition \ref{prop:fixrel}. 
Note that any birational map $f_{\delta,(a_i),(b_j)}$ is conjugate to $f_{\delta,(a_i/c),(b_j/c)}$ for any 
$c \in \mathbb{C}^*$ via the linear map $[x:y:z] \mapsto [c x : c y : z]$. 
Hence it is enough to show the existence of $(\delta,(a_i),(b_j))$ with 
$\sum_{i=1}^N 1/a_i - \sum_{j=1}^N 1/b_j \neq 0$ instead of $\sum_{i=1}^N 1/a_i - \sum_{j=1}^N 1/b_j = 1$. 
\par
For given real numbers $0=a_0 < a_1 < a_2 < \cdots < a_N$ and 
an $N$-tuple $b=(b_i) \in (\mathbb{R}^*)^N$ with $a_{i-1} < b_i < a_i$, put 
\[
g(x) = d \prod_{i=1}^N \Bigl(1 -\frac{x}{a_i} \Bigr), \qquad 
g_0(x) = \prod_{i=1}^N \Bigl(1- \frac{x}{b_i} \Bigr), 
\]
where $d=(1+\delta)^2/\delta \in [0,4]$ with $\delta \in S^1$. 
Moreover we assume that $0 < d <1$. 
Since $g(x)$ and $g_0(x)$ are polynomials of degree $n$ satisfying the relations 
$g(a_0)=d<1=g_0(a_0)$, $g(a_i)=0< (-1)^i g_0(a_i)$ and $(-1)^i g(b_i)<0=g_0(b_i)$ for $i \ge 1$, 
there is a unique real number $y_i \in (a_{i-1},b_i)$ such that $g(y_i)=g_0(y_i)$ for $i \in \{1,\dots,N\}$ 
(see Figure \ref{fig:fctn}). 
It is seen that $y_i=y_i(b)$ is continuous as a function of $b=(b_i) \in \prod_{i=1}^N(a_{i-1},a_i)$. 
\begin{figure}[t]
\begin{center}
\input{fctn.tex}
\end{center}
\caption{Two functions $g(x)$ and $g_0(x)$} 
\label{fig:fctn}
\end{figure}
\begin{lemma} \label{lem:fixrel}
Assume that $0 << d < 1$. Then there exists $b =(b_i) \in \prod_{i=1}^N(a_{i-1},a_i)$ such that 
$y_i(b)=x_i$ for any $i \in \{1,\dots,N \}$, where $x_i=(a_{i-1}+a_{i})/2$. 
Moreover, each component $b_i$ satisfies $\lim_{d \nearrow 1} b_i = a_i$. 
\end{lemma}
{\it Proof}. For $i \in \{1,\dots, N \}$ we put 
\[
s_i(x) = d \prod_{j=1}^i \Bigl( 1 - \frac{x}{a_j} \Bigr), \quad 
t_i(x) = \prod_{j=1}^i \Bigl( 1 - \frac{x}{a_j - \varepsilon_j} \Bigr), 
\]
where $\varepsilon_i$ is inductively determined by the relation $s_i(x_i)=t_i(x_i)$ (see also the following). 
We claim that $\varepsilon_i >0$ and 
$\varepsilon_i \searrow 0$ as $d \nearrow 1$. 
Indeed, if $i=1$, then the relation $s_1(x_1)=t_1(x_1)$ yields 
$\varepsilon_1=a_1 (1-d)(a_1 - x_1)/\{ a_1- d (a_1 - x_1) \} >0$, 
and $\varepsilon_1 \searrow 0$ as $d \nearrow 1$. 
Note that $d_2:=s_1(x_2) /t_1(x_2)$ satisfies $0 < d_2 < 1$ since $x_2 > a_1$, 
and $d_2 \nearrow 1$ as $d \nearrow 1$. 
Moreover for $i \ge 2$, assume that $d_i :=s_{i-1}(x_i)/t_{i-1}(x_i)$ 
satisfies $0 < d_i <1$, and $d_i \nearrow 1$ as $d \nearrow 1$. 
The relation 
\[
d_i \Bigl( 1- \frac{x_{i}}{a_{i}} \Bigr) = \frac{s_{i}(x_{i})}{t_{i-1}(x_{i})}
= \frac{t_{i}(x_{i})}{t_{i-1}(x_{i})} = \Bigl( 1- \frac{x_{i}}{a_{i} - \varepsilon_{i}} \Bigr) 
\]
yields $\varepsilon_{i} =a_i (1-d_i)(a_i - x_i)/\{ a_i- d_i (a_i - x_i) \} >0$ 
and $\varepsilon_i \searrow 0$ as $d \nearrow 1$. 
Similarly, $d_{i+1}= s_i(x_{i+1}) /t_i(x_{i+1})$ satisfies $0 < d_{i+1}<1$, and $d_{i+1} \nearrow 1$ as $d \nearrow 1$. 
Our claim is proved. 
\par
Now assume $0 << d <1$ so that $\varepsilon_i < a_i - x_i$. 
Regarding $y_i=y_i(b)$ as a function of $b=(b_i)$, 
we also claim that $y_i(b_1, \dots,b_{i-1},a_i - \varepsilon_i,b_{i+1}, \dots,b_N) < x_i$ 
for any $i \in \{1,\dots,N\}$ and any $(b_1, \dots,b_{i-1},b_{i+1}, \dots,b_N)$ with $(a_{j-1}<x_j <) a_j - \varepsilon_j < b_j < a_j$. 
Indeed, by putting  
\[
g_i(x) := \Bigl( 1 - \frac{x}{a_i - \varepsilon_i} \Bigr) \prod_{j \neq i}  \Bigl( 1 - \frac{x}{b_j} \Bigr),
\]
one has 
\[
\begin{array}{rl}
(-1)^{i-1} g_i(x_i) = & \displaystyle 
\Bigl( 1 - \frac{x_i}{a_i - \varepsilon_i} \Bigr) \prod_{j < i}  \Bigl( \frac{x_i}{b_j} -1 \Bigr) 
\prod_{j > i}  \Bigl( 1- \frac{x_i}{b_j}  \Bigr) \\[2mm]
< & \displaystyle
 \Bigl( 1 - \frac{x_i}{a_i - \varepsilon_i} \Bigr) \prod_{j < i}  \Bigl( \frac{x_i}{a_j - \varepsilon_j} -1 \Bigr) 
\prod_{j > i}  \Bigl( 1- \frac{x_i}{a_j}  \Bigr) \\[2mm]
= & \displaystyle 
d \Bigl( 1 - \frac{x_i}{a_i} \Bigr) \prod_{j < i}  \Bigl( \frac{x_i}{a_j} -1 \Bigr) 
\prod_{j > i}  \Bigl( 1- \frac{x_i}{a_j}  \Bigr) = (-1)^{i-1} g(x_i) 
\end{array}
\]
and $(-1)^{i-1} g_i(a_{i-1}) > 0 = (-1)^{i-1} g(a_{i-1})$, which yield the claim. 
\par
Finally we prove the existence of $b$ with $y_i(b)=x_i$. 
To this end, note that there is a root $z_i$ of $g(x)=g_0(x)$ such that $z_i \nearrow a_i$ as $b_i \nearrow a_i$.  
For $i=N$, the root $z_N$ must satisfy $z_N=y_N$ since $y_j \le a_{N-1}$ for $j \le N-1$. 
The above claim says that $y_N(b_1,\dots,b_{N-1},a_N -\varepsilon_N) < x_N$, 
which means that there is $b_N=b_N(b_1,\dots,b_{N-1}) \in (a_N-\varepsilon_N,a_N)$, depending continuously on $(b_j)_{j=1}^{N-1}$, 
such that $y_N(b_1,\dots,b_{N-1},b_N(b_1,\dots,b_{N-1})) = x_N$. 
Put $y_j(b_1,\dots,b_{N-1}) =y_j(b_1,\dots,b_{N-1},b_N(b_1,\dots,b_{N-1}))$, which is continuous with respect to $(b_j)_{j=1}^{N-1}$. 
Moreover for $i \le N-1$, we assume that $y_j=y_j(b_1,\dots,b_{i})$ satisfies $y_j=x_j$ for $j \ge i+1$. 
Similarly, $z_i$ must satisfy $z_i=y_i$ since $y_j \ge x_{i+1}$ for $j \ge i+1$ and $y_j \le a_{i-1}$ for $j \le i-1$. 
The above claim says that $y_i(b_1,\dots,b_{i-1},a_i -\varepsilon_i) < x_i$, 
which means that there is a continuous function 
$b_i=b_i(b_1,\dots,b_{i-1}) \in (a_i-\varepsilon_i,a_i)$ with $y_i(b_1,\dots,b_{i-1},b_i(b_1,\dots,b_{i-1})) = x_i$. 
Defining a continuous function $y_j(b_1,\dots,b_{i-1}) =y_j(b_1,\dots,b_{i-1},b_i(b_1,\dots,b_{i-1}))$, 
we can continue the induction. 
\par
To the end, there is $b=(b_i) \in \prod_{i=1}^N (a_i-\varepsilon_i,a_i)$ such that $y_i(b)=x_i$ for any $i \in \{1,\dots,N\}$. 
Since $\varepsilon_i \searrow 0$ as $d \nearrow 1$, we establish the lemma. 
\hfill $\Box$ \par\medskip 
\begin{lemma} \label{lem:siegel1} 
There exists $c_0 \in A$ such that the birational map $f$ determined by $c_0$ satisfies 
$\{ \mathrm{Tr}Df(w_l)\}^2/\mathrm{Det}Df(w_l) \in (0,4)$ for any 
$l \in \{1,\dots,N+2\}$. 
\end{lemma}
{\it Proof}.  
Under the notations mentioned in Lemma \ref{lem:fixrel},  
we can choose $0 <<d<1$ and $0 < b_1 < a_1 < b_2 < \cdots < b_N < a_N$ so that 
\[
\Bigl| \frac{1}{1-x_l/b_i} - \frac{1}{1-x_l/a_i} \Bigr| = \Bigl|\frac{(a_i-b_i) x_l}{(a_i-x_l)(b_i-x_l)} \Bigr| < \frac{1}{N} \quad 
(l \in \{1,\dots,N \}), \qquad 
1< \frac{a_i}{b_i}  < 2^{1/N}
\]
for any $i \in \{1,\dots,N\}$. 
Then from Lemma \ref{lem:fxtr1} and the fact $\mathrm{Det} Df(w_l)=\delta$, we have 
\[
\frac{\{\mathrm{Tr} Df (w_l)\}^2}{\mathrm{Det} Df(w_l)}= d \Bigl\{ 1 
+ \sum_{i=1}^N \Bigl( \frac{1}{1-x_l/b_i} -\frac{1}{1-x_l/a_i} \Bigr) \Bigr\}^2 \in (0,4) 
\]
for any $l \in \{1,\dots, N\}$. 
Choose $\delta \in S^1$ so that $d=(1+\delta)^2/\delta$. 
It follows from Lemma \ref{lem:fxtr2} and the fact $1< \beta_0/\alpha_0= \prod_{i=1}^N a_i/b_i<2$ that 
$\{\mathrm{Tr} Df (w_l)\}^2/\mathrm{Det} Df(w_l) \in [0,4]$ and then 
$\{\mathrm{Tr} Df (w_l)\}^2/\mathrm{Det} Df(w_l) \in (0,4)$ 
for $l \in \{N+1, N+2\}$ by slightly modifying the parameters if necessary. 
Thus, we have the desired parameters $c_0=(\delta,(a_i),(b_i))$.   
\hfill $\Box$ \par\medskip 
Next we consider the case $0 < b_0:=b_1 = \cdots = b_N < a_0:=a_1= \cdots =a_N$. 
Then the fixed points $w_l = [x_l:x_l:1]$ for $l \in \{1,\dots, N\}$ are given by the roots of 
$d (1-x_l/a_0)^N=(1-x_l/b_0)^N$ with $d= (1+\delta)^2/\delta \in [0,4]$, which yields 
\[
x_l = \frac{a_0 b_0 (1-\lambda_N \epsilon_N^l)}{a_0-b_0 \lambda_N \epsilon_N^l}, 
\]
where $\lambda_N:=d^{1/N} \ge 0$ and $\epsilon_N:=\cos (2 \pi/N) + i \sin (2 \pi/N)$ 
is a primitive $N$-th root of unity. 
Thus it follows from Lemma \ref{lem:fxtr1} and the fact $\mathrm{Det} Df(w_l)=\delta$ that 
\begin{equation} \label{eqn:ttd}
\begin{array}{l}
\displaystyle \frac{\{\mathrm{Tr} Df(w_l)\}^2}{\mathrm{Det} Df(w_l)} 
= d \Bigl\{ 1 - \frac{N}{a_0-b_0} 
(a_0 + b_0 - a_0  \lambda_N^{-1} \epsilon_N^{-l} - b_0  \lambda_N \epsilon_N^l) \Bigr\}^2 
\\~~~~ \displaystyle 
= d \Bigl\{ \Bigl( 1 - N \frac{a_0/b_0+1}{a_0/b_0-1} + N \frac{ \lambda_N^{-1} a_0 /b_0 + \lambda_N}{a_0/b_0-1} 
\cos \frac{2 \pi l}{N} \Bigr)- i N \frac{\lambda_N^{-1} a_0 /b_0 - \lambda_N}{a_0/b_0-1} \sin \frac{2 \pi l}{N}  \Bigr\}^2 
\end{array}
\end{equation}
for $l \in \{1,\dots, N\}$. Moreover from Lemma \ref{lem:fxtr2}, one has 
\[
\frac{\{\mathrm{Tr} Df(w_l)\}^2}{\mathrm{Det} Df(w_l)} \notin [0,4] 
\Longleftrightarrow \Bigl( \frac{a_0}{b_0} \Bigr)^N \notin [0,4]
\]
for $l \in \{ N+1,N+2 \}$. 
\begin{lemma} \label{lem:siegel2} 
There exists $c_* \in A$ such that the birational map $f$ determined by $c_*$ satisfies 
$\{ \mathrm{Tr}Df(w_l)\}^2/\mathrm{Det}Df(w_l) \notin [0,4]$ for any 
$l \in \{1,\dots,N+2\}$. 
\end{lemma}
{\it Proof}. 
First we assume that $d=1/4^2$ and $a_0/b_0=4^{1/N}$ in the above notations. 
If $l \in \{1,\dots,N\}$, the only possibilities for $\{\mathrm{Tr} Df(w_l)\}^2/\mathrm{Det} Df(w_l)$ 
to become a nonnegative real number occur when $(\cos 2 \pi l/N, \sin 2 \pi l/N) = (\pm 1,0)$ in (\ref{eqn:ttd}). 
On the other hand, in the case $(\cos 2 \pi l/N, \sin 2 \pi l/N) = (\pm 1,0)$, it is seen that 
$\{\mathrm{Tr} Df(w_l)\}^2/\mathrm{Det} Df(w_l)>4$. 
Indeed, when $(\cos 2 \pi l/N, \sin 2 \pi l/N) = (1,0)$, one has  
\[
\frac{\{\mathrm{Tr} Df(w_l)\}^2}{\mathrm{Det} Df(w_l)} 
=  \frac{1}{4^2} \Bigl( 1 - N \frac{4^{1/N}+1}{4^{1/N}-1} + N \frac{ 4^{3/N} + 4^{-2/N}}{4^{1/N}-1} \Bigr)^2
= \Bigl( 2 +\frac{N}{4} g(N) \Bigr)^2,
\]
where $g(N):= (4^{2/N}-4^{-2/N})+(4^{1/N}-4^{-1/N})-7/N$. 
The function $g(N)$ satisfies $g(N)>0$ for any $N \ge 1$, 
as $g(N)$ is monotone decreasing in $N$ and $\lim_{N \to \infty} g(N)=0$. 
Thus we have $\{\mathrm{Tr} Df(w_l)\}^2/\mathrm{Det} Df(w_l)>4$. 
The case $(\cos 2 \pi l/N, \sin 2 \pi l/N) = (- 1,0)$ can be treated in a similar manner. 
Thus the condition $\{\mathrm{Tr} Df(w_l)\}^2/\mathrm{Det} Df(w_l) \notin [0,4]$ holds 
for any $l \in \{1,\dots,N\}$. 
\par
Now since $\{\mathrm{Tr} Df(w_l)\}^2/\mathrm{Det} Df(w_l)$ continuously depends on the parameters $(\delta,a,b) \in A$, 
with the above condition, we slightly modify the parameters so that 
$0 < b_1 < \cdots < b_N < b_0 <a_0 < a_1 < \cdots <a_N$, 
which means that $\beta_0/\alpha_0>(a_0/b_0)^N=4$ and thus 
$\{\mathrm{Tr} Df(w_l)\}^2/\mathrm{Det} Df(w_l) \notin [0,4]$ for any $l \in \{1,\dots,N+2\}$. 
By fixing $\delta \in S^1$ with $d= (1+\delta)^2/\delta$, 
we show the existence of $c_*=(\delta, a,b) \in A$.  
\hfill $\Box$ \\[2mm]
{\it Proof of Proposition \ref{prop:fixrel}}. 
Note that $\{\mathrm{Tr} Df(w_l)\}^2/\mathrm{Det} Df(w_l)$ 
continuously depends on the parameters $(\delta,a,b) \in A$.  
Hence the proposition is the consequence of Lemmas \ref{lem:siegel1} and \ref{lem:siegel2}. 
\hfill $\Box$ \par\medskip 

\end{document}